\documentclass[11pt]{article}
  \usepackage{graphicx,wrapfig}
 \usepackage{flafter}
 \usepackage{leftidx}
 \usepackage{mathrsfs}
 \usepackage{amssymb,amsmath}
 \usepackage{bbding}
\usepackage{fancyhdr}%
\usepackage{mathrsfs}
\usepackage{color}
\usepackage{stmaryrd}
\usepackage{amsfonts}
\usepackage{latexsym}
\usepackage{psfrag}
\usepackage{graphicx,subfigure}
\usepackage{fancyhdr,graphicx}
\usepackage{multicol}
\usepackage{dsfont}
\usepackage{bbm}
\usepackage{booktabs}
\usepackage[center]{caption2}
\usepackage{cite}
\usepackage{epstopdf}

\allowdisplaybreaks

\renewcommand{\leq}{\leqslant}
\renewcommand{\theequation}{\thesubsection.\arabic{equation}}

\date{}
\def\theequation{\arabic{section}.\arabic{equation}}
\newtheorem{theorem}{Theorem}[section]
\newtheorem{lemma}[theorem]{Lemma}
\newtheorem{remark}[theorem]{Remark}

\newcommand{\zd}{\,\mathrm{d}}
\newcommand{\myinner}[1]{\left\langle#1\right\rangle}
\newcommand{\myinnerb}[1]{\big\langle#1\big\rangle}
\newcommand{\mynorm}[1]{\left\|#1\right\|}
\newcommand{\mynormb}[1]{\big\|#1\big\|}
\newcommand{\mynormt}[1]{\|#1\|}
\newcommand{\abs}[1]{\left|#1\right|}
\newcommand{\absb}[1]{\big|#1\big|}
\newcommand{\absB}[1]{\Big|#1\Big|}

\newcommand{\bra}[1]{\left(#1\right)}
\newcommand{\brab}[1]{\big(#1\big)}
\newcommand{\braB}[1]{\Big(#1\Big)}
\newcommand{\brat}[1]{(#1)}
\newcommand{\kbra}[1]{\left[#1\right]}
\newcommand{\kbrab}[1]{\big[#1\big]}
\newcommand{\kbraB}[1]{\Big[#1\Big]}

\newcommand{\fd}[1]{\mathcal{D}^{#1}_t}
\newcommand{\dfd}[1]{D^{#1}_\tau}
\newcommand{\ffd}[1]{D^{#1}_f}
\newcommand{\diff}{\nabla_\tau}
\newcommand{\Ass}[1]{\textbf{\upshape Ass#1}}
\newcommand{\vecx}{\boldsymbol{x}}
\newcommand{\taumax}{\tau}
\newcommand{\rhomax}{\rho}
\newcommand{\defeq}{:=}
\newcommand{\Err}[1]{\widetilde{\Pi_{#1}}}

\begin{document}

\title{Unconditional convergence of a fast two-level linearized algorithm
for semi\/linear subdiffusion equations}
\author{Hong-lin Liao\thanks{
Department of Mathematics, Nanjing University of Aeronautics and Astronautics,
Nanjing, 211106, P. R. China.
E-mail: liaohl@csrc.ac.cn. He is partially supported by NSFC grant 11372354
and a grant DRA2015518 from 333 High-level Personal Training Project of Jiangsu Province.}
\quad Yonggui Yan\thanks{
Beijing Computational Science Research Center (CSRC), Beijing, 100193, P. R. China.
 E-mail: yan\_yonggui@csrc.ac.cn. }
\quad Jiwei Zhang\thanks{
Beijing Computational Science Research Center (CSRC),
Beijing, 100193, P. R. China.  E-mail: jwzhang@csrc.ac.cn.
He is partially supported by NSFC under grants 11771035, 91430216 and NSAF U1530401.}
}

 \maketitle\normalsize

\begin{abstract}

A fast two-level linearized scheme with unequal time-steps is constructed and analyzed
for an initial-boundary-value problem of semilinear subdiffusion equations.
The two-level fast L1 formula of the Caputo derivative is derived based on
the sum-of-exponentials technique. The resulting fast algorithm
is computationally efficient in long-time simulations because it
significantly reduces the computational cost $O(MN^2)$ and storage $O(MN)$
for the standard L1 formula to $O(MN\log N)$ and $O(M\log N)$,
respectively, for $M$ grid points in space and $N$ levels in time.
The nonuniform time mesh would be graded to handle the typical singularity
of the solution near the time $t=0$,
and Newton linearization is used to approximate the nonlinearity term.
Our analysis relies on three tools:
a new discrete fractional Gr\"{o}nwall inequality, a global
consistency analysis and a discrete $H^2$ energy method.
A sharp error estimate reflecting the regularity of solution is established without
any restriction on the relative diameters of the temporal and spatial mesh sizes.
Numerical examples are provided to demonstrate the effectiveness of our approach
and the sharpness of error analysis.
\vspace{1em}\\
\emph{Keywords}:\;\;semilinear subdiffusion equation; two-level L1 formula;
discrete fractional Gr\"{o}nwall inequality;
discrete $H^2$ energy method; unconditional convergence

\end{abstract}


\section{Introduction}
\setcounter{equation}{0}

A two-level linearized method is considered to numerically solve the following semilinear subdiffusion
equation on a bounded domain
\begin{subequations}\label{prob}
\begin{align}
\fd{\alpha}u &= \Delta u+f(u)\quad\text{for $\vecx\in\Omega$ and $0<t\leq T$},\label{Problem-1}\\
u&=u^0(\vecx) \quad\text{for $\vecx \in\Omega$ and $t=0$},\label{Problem-2}\\
u&= 0\quad\text{for $\vecx\in\partial\Omega$ and $0<t\leq T$},\label{Problem-3}
\end{align}
\end{subequations}
where $\partial\Omega$ is the boundary of $\Omega:=(x_{l},x_{r})\times (y_{l},y_{r})$,
and the nonlinear function $f(u)$ is smooth.
In~\eqref{Problem-1} $\fd{\alpha}={}_{~0}^{C}\!\mathcal{D}_t^{\alpha}$ denotes
the  Caputo  fractional derivative of order $\alpha$:   
\begin{align}\label{CaputoDef}
\brat{\fd{\alpha}v}(t)
\defeq\int_{0}^{t} \omega_{1-\alpha} (t-s) v'(s)\zd{s},\quad
0<\alpha<1,
\end{align}
where the weakly singular kernel $\omega_{1-\alpha} (t-s)$ is defined by
$\omega_{\mu}(t):=t^{\mu-1}/{\Gamma(\mu)}$. It is easy to verify
$\omega_{\mu}'(t)=\omega_{\mu-1}(t)$ and $\int_0^t\omega_{\mu}(s)\zd s=\omega_{\mu+1}(t)$ for $t>0$.

In any numerical methods for solving nonlinear fractional diffusion equations \eqref{Problem-1},
a key consideration is the singularity of the solution near the time $t=0$, see \cite{JinLiZhou:2017,LiYiChen:2016,ZhangSunLiao:2014,MustaphaMustapha:2010}.
For example, under the assumption that the nonlinear function $f$ is Lipschitz continuous
and the initial data $u^0\in H^2(\Omega)\cap H_0^1(\Omega)$,
Jin et al. \cite[Theorem 3.1]{JinLiZhou:2017} prove that problem \eqref{prob} has a unique solution $u$ for which
$u\in C\bra{[0,T];H^2(\Omega)\cap H_0^1(\Omega)}$, $\fd{\alpha}u\in C\bra{[0,T];L^2(\Omega)}$ and $
\partial_t u \in L^2(\Omega)$  with
$\|\partial_t u(t)\|_{L^2(\Omega)} \le C_ut^{\alpha-1}$ for $0<t\leq T$,
where $C_u>0$ is a constant independent of $t$ but may depend on~$T$.
Their analysis of numerical methods for solving~\eqref{prob} is applicable to both the L1 scheme
and backward Euler convolution quadrature on a uniform time grid of diameter~$\tau$;
 a lagging linearized technique is used to handle the nonlinearity $f(u)$, and \cite[Theorem 4.5]{JinLiZhou:2017}
 shows that the discrete solution  is $O(\tau^{\alpha})$ convergent in $L^\infty(L^2(\Omega))$.

This work may be considered as a continuation of \cite{LiaoLiZhangZhao:2017}, in which a sharp error estimate
for the L1 formula on nonuniform meshes was obtained for linear subdiffusion-reaction equations based on
a discrete fractional Gr\"{o}nwall inequality and a global consistency analysis. In this paper,
we combine the L1 formula and the sum-of-exponentials (SOEs) technique to develop a
one-step fast difference algorithm for the nonlinear subdiffusion problem \eqref{prob} by
using the Newton's linearization to approximate nonlinear term,
and present the corresponding sharp error estimate of the proposed scheme
without any restriction on the relative diameters of temporal and spatial mesh sizes.

It is known that the Caputo fractional derivative involves a convolution kernel.
The total number of operations required to evaluate the sum of L1 formula is proportional
to $O(N^2)$, and the active memory to $O(N)$ with $N$ representing the total time steps,
which is prohibitively expensive for the practically large-scale and long-time simulations.
Recently, a simple fast algorithm based on SOEs approximation is proposed to significantly
reduce the computational complexity to $O(N \log N)$ and $O(\log N)$ when the final time $T\gg 1$,
see \cite{Li:2010, JiangZhangZhangZhang:2017}.  Another fast algorithm for the evaluation of
the fractional derivative has been proposed in \cite{BH17}, where the compression is carried
out in the Laplace domain by solving the equivalent ODE with some one-step A-stable scheme.
In this paper,  we develop a fast two-level L1 formula
by combining a nonuniform mesh suited to the initial singularity
with a fast time-stepping algorithm for the historical memory in~\eqref{CaputoDef}.
This scheme would be also useful to develop efficient parallel-in-time  algorithms
for time-fractional differential equations  \cite{XuHesthavenChen:2015}.

On the other hand, the nonlinearity of the problem also results in the difficulty for the numerical analysis.
To establish an error estimate of the two-level linearized scheme at time~$t_n$, it requires to prove
the boundedness of the numerical solution at the previous time levels via
$\mynormt{u^{n-1}}_{\infty}\leq C_u$.
Traditionally it is done using mathematical induction and some inverse estimate, namely,
\begin{align*}
\mynormt{u^{n-1}}_{\infty}\leq&\, 
\mynormt{U^{n-1}}_{\infty}+h^{-1}\mynormt{U^{n-1}-u^{n-1}}
\leq
\mynormt{U^{n-1}}_{\infty}+C_uh^{-1}\brab{\tau^{\beta}+h^2}.
\end{align*}
This leads to that a time-space grid
restriction $\tau=O(h^{1/\beta})$ is required in the theoretical analysis even though
it is nonphysical and may be unnecessary in numerical simulations. In this paper, we will
extend the discrete $H^2$ method developed in \cite{LiaoSun:2010,LiaoSunShi:2010,LiaoSunShi:2010B}
to prove unconditional convergence of our fully discrete solution without the restriction conditions of
between mesh sizes $\tau$ and~$h$ comparing with the traditional method.  The main idea of
discrete $H^2$ energy method is to separately treat the temporal and spatial truncation errors.
This simple implementation avoids some nonphysical time-space grid restrictions in the error analysis.
A related approach in a finite element setting are discussed in~\cite{LiSun:2013,LiSun:2013sinum,LiWangSun:2014}.

The convergence rate of L1 formula for the Caputo derivative is limited by the smoothness of the solution.
The analysis here is based on the following assumptions on the solution
\begin{align}\label{SolutionRegularityAssumption}
\mynorm{u}_{H^{4}(\Omega)}\leq C_u,\;
\mynorm{\partial_{t}u}_{H^{4}(\Omega)}\leq C_u(1+t^{\sigma-1}) \text{ and }
\mynorm{\partial_{tt}u}_{H^{2}(\Omega)}\leq C_u(1+t^{\sigma-2})
\end{align}
for $0< t \le T$, where $\sigma\in(0,1)\cup(1,2)$ is a regularity parameter.
To resolve the singularity at~$t=0$, it is reasonable to use a nonuniform mesh that concentrates  grid points near $t=0$,
see \cite{Brunner:book,BrunnerLingYamamoto:2010,LiaoLiZhangZhao:2017,StynesRiordanGracia:2017}.
We make the following assumption on the time mesh:
\begin{description}
\item[\Ass{G}.] Let $\gamma\ge1$ be a user-chosen parameter. There is a constant $C_{\gamma}>0$, independent of $k$,  such that
$\tau_k\le C_{\gamma}\tau\min\{1,t_k^{1-1/\gamma}\}$
for~$1\le k\le N$ and $t_{k}\leq C_{\gamma}t_{k-1}$ for~$2\le k\le N$.
\end{description}
Since $\tau_1=t_1$, \Ass{G} implies that $\tau_1=O(\tau^{\gamma})$, while for those $t_k$ bounded
away from $t=0$ one has $\tau_k=O(\tau)$. The parameter~$\gamma$ controls the extent to which the
grid points are concentrated near~$t=0$: increasing $\gamma$ will decrease the time-step sizes
near $t=0$ and so move mesh points closer to $t=0$.
A simple example of a family of meshes satisfying \Ass{G} is the graded grid $t_k=T(k/N)^{\gamma}$,
which is discussed in~\cite{Brunner:book,LiaoLiZhangZhao:2017,StynesRiordanGracia:2017}.  Although nonuniform meshes are flexible and reasonably convenient for practical
implementation, they can significantly complicate the numerical analysis of
schemes, both with respect to stability and consistency. In this paper, our analysis will rely on a generalized fractional Gr\"{o}nwall inequality~\cite{LiaoMcLeanZhangZhao:2018},
which would be applicable for any discrete fractional derivatives having the discrete convolution form.

Throughout the paper, any subscripted $C$, such as $C_u$, $C_{\gamma}$, $C_{\Omega}$,
 $C_v$, $C_0$ and $C_{F}$, denotes a generic positive constant, not necessarily
 the same at different occurrences, which is always dependent on the given data
 and the solution but independent of the time-space grid steps.
 The  paper is organized as follows. Section~\ref{sec:two-level} presents the
two-level fast L1 formula and the corresponding linearized fast scheme.
The global consistency analysis of fast L1 formula and the Newton's linearization
is presented in Section~\ref{sec:consistency}. A sharp error estimate
for the linearized fast scheme is proved in Section~\ref{sec:error}.
Two numerical examples in Section~\ref{sec:numerical} are given to demonstrate the sharpness of our analysis.

\section{A two-level fast method}\label{sec:two-level}
\setcounter{equation}{0}

We approximate the Caputo fractional  derivative \eqref{CaputoDef} on a (possibly nonuniform) time mesh
$0=t_0<\cdots<t_{k-1}<t_{k}<\cdots<t_N=T$, with the time-step sizes
$\tau_k\defeq t_k-t_{k-1}$ for $1\le k\le N$, the maximum time-step $\taumax=\max_{1\leq k\leq N}\tau_k$
and the step size ratios~$\rho_k\defeq\tau_k/\tau_{k+1}$ for $1\le k\le N-1$.
In space we use a standard finite difference method on a tensor product grid.
 Let $M_1$ and $M_2$ be two positive integers. Set $h_{1}=(x_{r}-x_{l})/M_1, \ h_{2}=(y_{r}-y_{l})/ M_2$
 and the maximum spatial length $h=\max\{h_1,h_2\}$. Then the fully discrete spatial grid
 $\bar{\Omega}_h :=\{\vecx _{h}=(x_{l} + ih_1, y_{l} + jh_2)\,|\, 0\leq i\leq M_1, 0\leq j\leq M_2\}$. Set ${\Omega}_h=\bar\Omega_h\cap\Omega$ and the boundary $\partial
\Omega_h=\bar\Omega_h\cap\partial\Omega$. Given a grid function $v = \{v_{ij}\}$, define
\[
v_{i-\frac12,j}=\bra{v_{i,j}+v_{i-1,j}}/{2}, \ \delta_xv_{i-\frac12,j}=\bra{v_{i,j}-v_{i-1,j}}/{h_1}, \ \delta_x^2v_{ij}=\brab{\delta_xv_{i+\frac12,j}-\delta_xv_{i-\frac12,j}}/{h_1}.
\]
Difference operators $v_{i,j-\frac12}, \ \delta_yv_{i,j-\frac12}, \ \delta_x\delta_yv_{i-\frac12,j-\frac12}$
and $\delta_y^2v_{ij}$ can be defined analogously.
The second-order approximation of $\Delta v(\vecx _{h})$ for $\vecx _h\in\Omega_h$
is $\Delta_hv_{h} := (\delta_x^2+\delta_y^2)v_{h}$.
Let $\mathcal{V}_{h}$ be the space of grid functions,
$\mathcal{V}_{h}=\big\{v=\bra{v_h}_{\vecx _h\in\bar{\Omega}_h}\,\big|\,
v_h=0\;\mbox{for}\;\vecx_h\in\partial\Omega_h\big\}$.
 For $v, w \in \mathcal{V}_h$, define the discrete inner product
 $\myinner{v,w}=h_1h_2\sum_{\vecx_h\in\Omega_h}v_{h}w_{h}$,
 the $L^2$ norm $\|v\|=\sqrt{\myinner{v,v}}$,
 the $H^1$ seminorm $\|\nabla_h v\|=\sqrt{\mynormt{\delta_xv}^2+\mynormt{\delta_yv}^2}$ and the maximum norm
$\|v\|_{\infty}=\max_{\vecx_h\in\Omega_h}\abs{v_h}$. For any $v\in\mathcal{V}_{h}$,
by~\cite[Lemmas 2.1, 2.2 and 2.5]{LiaoSun:2010} there exists a constant $C_{\Omega}>0$ such that
\begin{align}\label{discreteNormsEmbedding}
\|v\|\leq C_{\Omega}\|\nabla_hv\|,
\quad \|\nabla_hv\|\leq C_{\Omega}\|\Delta_hv\|,
\quad \|v\|_{\infty}\leq C_{\Omega}\|\Delta_hv\|.
 \end{align}

\subsection{A fast variant of the L1 formula}

On our nonuniform mesh, the standard L1 approximation of the Caputo derivative is
\begin{align}\label{L1formulaNonuniform}
(\dfd{\alpha}v)^n\defeq\sum_{k=1}^n\frac{1}{\tau_k}\int_{t_{k-1}}^{t_{k}}\omega_{1-\alpha}(t_n-s)\nabla_{\tau} v^k\zd s
= \sum_{k=1}^na_{n-k}^{(n)}\nabla_{\tau}v^k\,,
\end{align}
where $\nabla_{\tau}v^k:=v^k-v^{k-1}$ and the convolution kernel $a_{n-k}^{(n)}$ is defined by
\begin{align}\label{Coefficient1-L1formulaNonuniform}
a_{n-k}^{(n)}\defeq\frac{1}{\tau_k}\int_{t_{k-1}}^{t_{k}}\!\!\omega_{1-\alpha}(t_n-s)\zd s
=\frac{1}{\tau_k}\kbra{\omega_{2-\alpha}(t_n-t_{k-1})-\omega_{2-\alpha}(t_n-t_{k})},\quad 1\leq k\leq n.
\end{align}

\begin{lemma}\label{lemma: L1 kernel}
For fixed integer $n\geq2$, the convolution kernel $a^{(n)}_{n-k}$ of~\eqref{Coefficient1-L1formulaNonuniform} satisfies
\begin{itemize}
  \item[(i)] $\displaystyle a_{n-k-1}^{(n)}>\omega_{1-\alpha}(t_n-t_{k})>a_{n-k}^{(n)},\quad 1\leq k\leq n-1;$
  \item[(ii)]$\displaystyle a^{(n)}_{n-k-1}-a^{(n)}_{n-k}>\tfrac12\kbra{\omega_{1-\alpha}(t_{n}-t_k)
  -\omega_{1-\alpha}(t_{n}-t_{k-1})},\quad 1\leq k\leq n-1.$
\end{itemize}
\end{lemma}
\begin{proof} The integral mean-value theorem yields (i) directly;
see~\cite{ZhangSunLiao:2014,LiaoLiZhangZhao:2017}.
For any function $q\in C^2[t_{k-1},t_{k}]$, let $\Pi_{1,k}q$ be the linear interpolant
of $q(t)$ at  $t_{k-1}$ and $t_k$. Let
$\Err{1,k}q := q - \Pi_{1,k}q$ be the error in this interpolant.
For $q(s)=\omega_{1-\alpha}(t_{n}-s)$ one has $q''(s)=\omega_{-\alpha-1}(t_{n}-s)>0$ for $0<s<t_n$, so
the Peano representation of the interpolation error \cite[Lemma 3.1]{LiaoLiZhangZhao:2017}  shows that
$\int_{t_{k-1}}^{t_{k}}\brab{\Err{1,k}q}(s)\zd s<0$.
Thus the definition \eqref{Coefficient1-L1formulaNonuniform} of $a^{(n)}_{n-k}$ yields
\begin{align*}
a^{(n)}_{n-k}-\frac12\omega_{1-\alpha}(t_{n}-t_k)
    -\frac12\omega_{1-\alpha}(t_{n}-t_{k-1})
    =\frac1{\tau_k}\int_{t_{k-1}}^{t_{k}}\brab{\Err{1,k}q}(s)\zd s<0,\quad 1\leq k\leq n-1.
\end{align*}
Subtract this inequality from (i) to obtain (ii) immediately.
\end{proof}

As the L1 formula \eqref{L1formulaNonuniform} involves the solution at all previous time-levels,
it is computationally inefficient to directly evaluate it  when solving
the fractional diffusion problem \eqref{prob} using time-stepping.
We therefore use the SOEs approach
of~\cite{JiangZhangZhangZhang:2017,Li:2010,YanSunZhang:2017} to develop
a fast L1 formula. A basic result of the SOE approximation
(see  \cite[Theorem 2.5]{JiangZhangZhangZhang:2017}
or \cite[Lemma 2.2]{YanSunZhang:2017}) is the following:

\begin{lemma}\label{lemma:sum-of-exponentials}
Given $\alpha\in(0,1)$, an absolute tolerance error $\epsilon\ll1$,
a cut-off time $\Delta{t}>0$ and a final time $T$, there exists a positive
integer $N_{q}$, positive quadrature nodes $\theta^{\ell}$ and positive
weights
$\varpi^{\ell}$ $(1\leq \ell\leq N_q)$ such that
\[
\absB{\omega_{1-\alpha}(t)-\sum_{\ell=1}^{N_q}\varpi^{\ell}e^{-\theta^{\ell}t}}
\leq \epsilon \quad \forall\, t\in[\Delta{t},T],
\]
where the number $N_q$  of quadrature nodes satisfies
$$N_q=O\bra{\log\frac1{\epsilon}\braB{\log\log\frac1{\epsilon}+\log\frac{T}{\Delta{t}}}
+\log\frac1{\Delta{t}}\braB{\log\log\frac1{\epsilon}+\log\frac{1}{\Delta{t}}}}.$$
\end{lemma}

After that, we divide the fractional Caputo  derivative $\brat{\fd{\alpha}v}(t_n)$
of~\eqref{CaputoDef} into a sum of a local part (an integral over $[t_{n-1}, t_n]$) and
a history part (an integral over $[0, t_{n-1}]$), then approximate $v'$ by linear interpolation in the local part (similar to the standard L1 method) and use
the SOE technique of Lemma~\ref{lemma:sum-of-exponentials} to approximate the kernel
$\omega_{1-\alpha} (t-s)$ in the history part.  It yields
\begin{align*}
\brab{\fd{\alpha}u}(t_n)
\approx&\, \int_{t_{n-1}}^{t_{n}}\omega_{1-\alpha}(t_n-s)\frac{\nabla_{\tau}u^n}{\tau_n}\zd s
+\int_0^{t_{n-1}}\sum_{\ell=1}^{N_q}\varpi^{\ell}e^{-\theta^{\ell}(t_n-s)}u'(s)\zd s\\
=&\,a_{0}^{(n)}\nabla_{\tau}u^n+\sum_{\ell=1}^{N_q}\varpi^{\ell}e^{-\theta^{\ell}\tau_n}\mathcal{H}^{\ell}(t_{n-1}),\quad n\geq1,
\end{align*}
where $\mathcal{H}^{\ell}(t_k) := \int_0^{t_{k}}e^{-\theta^{\ell}(t_k-s)}u'(s)\zd s$
with $\mathcal{H}^{\ell}(t_0)=0$ for $1\leq \ell\leq N_q\,$.
To compute $\mathcal{H}^{\ell}(t_k)$ efficiently we apply linear interpolation
in each cell $[t_{k-1}, t_{k}]$, obtaining
\begin{align*}
\mathcal{H}^{\ell}(t_{k})=&\,e^{-\theta^{\ell}\tau_{k}}\mathcal{H}^{\ell}(t_{k-1})
+\int_{t_{k-1}}^{t_{k}}e^{-\theta^{\ell}(t_{k}-s)}u'(s)\zd s
\approx e^{-\theta^{\ell}\tau_{k}}\mathcal{H}^{\ell}(t_{k-1})
+b^{(k,\ell)}\nabla_{\tau}u^{k},
\end{align*}
where the positive coefficient is given by
\begin{align}\label{Coefficient1-fastL1formulaNonuniform}
b^{(k,\ell)}\defeq\frac{1}{\tau_k}\int_{t_{k-1}}^{t_{k}}e^{-\theta^{\ell}(t_{k}-s)}\zd s,
\quad k\geq1, \;1\leq \ell\leq N_q\,.
\end{align}

In summary, we now have the two-level fast L1 formula
\begin{subequations}\label{fastL1}
\begin{align}\label{fastL1formulaNonuniform}
 \brat{\ffd{\alpha}u}^n\defeq
  a_{0}^{(n)}\nabla_{\tau}u^n+\sum_{\ell=1}^{N_q}\varpi^{\ell}e^{-\theta^{\ell}\tau_n}H^{\ell}(t_{n-1}), \quad n\geq1,
\intertext{where $H^{\ell}(t_{k})$ satisfies $H^{\ell}(t_{0})=0$ and the recurrence relationship}
\label{fastL1formulaRecurrence}
H^{\ell}(t_{k})=e^{-\theta^{\ell}\tau_{k}}H^{\ell}(t_{k-1})
+b^{(k,\ell)}\nabla_{\tau}u^{k},\quad k\geq1, \;1\leq \ell\leq N_q\,.
\end{align}
\end{subequations}

\subsection{The two-level linearized scheme}

Write $U_{h}^n=u(\vecx_h,t_n)$ for $\vecx_h\in\bar{\Omega}_h$, $0\leq n\leq N$.
Let $u_{h}^n$ be the discrete approximation of~$U_{h}^n$.
Using the fast L1 formula \eqref{fastL1} and Newton linearization, we obtain a
linearized scheme for the problem \eqref{prob}: find  $\{u_h^{N}\}\in\mathcal{V}_{h}$ such that
\begin{subequations}\label{diffScheme}
\begin{align}
\brat{\ffd{\alpha}u_h}^n=&\,\Delta_hu_h^n+f(u_h^{n-1})
+f'(u_h^{n-1})\nabla_{\tau}u_h^{n}\,,\quad
\vecx_h\in\Omega_h, \; 1\leq n\leq N;\label{differenceScheme-1}\\
u_h^0=&\,u^0(\vecx_h),\quad \vecx_h\in\bar{\Omega}_h\,.\label{differenceScheme-2}
\end{align}
\end{subequations}

Note that, the Newton linearization of a general nonlinear function $f=f(\vecx,t,u)$ at $t=t_n$ takes
the form $f(\vecx_h,t_n,u_h^n)\approx f(\vecx_h,t_n,u_h^{n-1}) +f_u'(\vecx_h,t_n,u_h^{n-1})\nabla_{\tau}u_h^{n}\,.$
The scheme \eqref{diffScheme} is a two-level procedure for computing $\{u_h^{n}\}$, since
 \eqref{differenceScheme-1} can be reformulated as
\begin{align}
\kbra{a_{0}^{(n)}-\Delta_h-f'(u_h^{n-1})}\nabla_{\tau}u_h^{n}
=&\,\Delta_hu_h^{n-1}+f(u_h^{n-1})
-\sum_{\ell=1}^{N_q}\varpi^{\ell}e^{-\theta^{\ell}\tau_n}H_h^{\ell}(t_{n-1}),\label{computProcess-1}\\
H_h^{\ell}(t_{n})=&\,e^{-\theta^{\ell}\tau_{n}}H_h^{\ell}(t_{n-1})
+b^{(n,\ell)}\nabla_{\tau}u_h^{n}\,,\quad  1\leq \ell\leq N_q.\label{computProcess-2}
\end{align}
Thus, once the solution $\{u_h^{n-1},\;H_h^{\ell}(t_{n-1})\}$ at the previous time-level $t_{n-1}$ is available,
the current solution $\{u_h^{n}\}$ can be found by \eqref{computProcess-1} with a fast matrix solver
and the historic term $\{H_h^{\ell}(t_{n})\}$ will be updated explicitly by the recurrence formula \eqref{computProcess-2}.

\begin{remark}\label{rem:efficiency}
At each time level the scheme \eqref{diffScheme} requires $O(MN_q)$ storage and $O(MN_q)$ operations,
where $M=M_1M_2$ is the total number of spatial grid points. Given a tolerance error $\epsilon=\epsilon_0$,
by virtue of Lemma~\ref{lemma:sum-of-exponentials}, the number of quadrature nodes $N_q=O(\log N)$
if the final time $T\gg1$. Hence our new method is computationally efficient since it computes
the final solution using in total $O(M\log N)$ storage and $O(MN\log N)$ operations.
\end{remark}

\subsection{Discrete fractional Gr\"{o}nwall inequality}

Our analysis relies on a generalized discrete fractional Gr\"{o}nwall inequality~\cite{LiaoMcLeanZhangZhao:2018},
which is applicable for any discrete fractional derivative having the discrete convolution form~
\begin{align}\label{DiscreteCaputoDef}
\brat{\fd{\alpha}v}^n
\approx\sum_{k=1}^n A^{(n)}_{n-k}(v^k-v^{k-1}),\quad
1\leq n\leq N,
\end{align}
provided that $A^{(n)}_{n-k}$ and the time-steps $\tau_n$ satisfy the following three assumptions:
\begin{description}
\item[\Ass{1}.] The discrete kernel is monotone, that is, $A^{(n)}_{k-2}\ge A^{(n)}_{k-1}>0$ for~$2\leq k\leq n\leq N$.
\item[\Ass{2}.] There is a constant~$\pi_A>0$ such that $A^{(n)}_{n-k}\ge\frac{1}{\pi_A}\int_{t_{k-1}}^{t_k}
\frac{\omega_{1-\alpha}(t_n-s)}{\tau_k}\zd{s}$ for $1\le k\le n\leq N$.
\item[\Ass{3}.] There is a constant~$\rhomax>0$ such that the time-step ratios $\rho_k\le\rhomax$
for~$1\le k\le N-1$.
\end{description}

The  complementary discrete kernel~$P^{(n)}_{n-k}$ was introduced by
Liao et al.~\cite{LiaoLiZhangZhao:2017,LiaoMcLeanZhangZhao:2018};
it satisfies the following identity
\begin{equation}\label{eq: P A}
\sum_{j=k}^nP^{(n)}_{n-j}A^{(j)}_{j-k}\equiv1\quad\text{for $1\le k\le n\le N$.}
\end{equation}
Rearranging this identity yields a recursive formula that defines $P^{(n)}_{n-k}$ :
\begin{align}\label{discreteConvolutionKernel-RL}
P_{0}^{(n)}\defeq1/{A_0^{(n)}},\quad
P_{n-j}^{(n)}\defeq1/{A_0^{(j)}}
\sum_{k=j+1}^{n}\braB{A_{k-j-1}^{(k)}-A_{k-j}^{(k)}}P_{n-k}^{(n)}\,,\quad 1\leq j\leq n-1.
\end{align}
From \cite[Lemma~2.2]{LiaoMcLeanZhangZhao:2018} we see that~$P^{(n)}_{n-k}$ is well-defined
and non-negative if the assumption \Ass{1} holds true. Furthermore, if \Ass{2} holds true, then
\begin{equation}\label{eq: P bound}
\sum_{j=1}^nP^{(n)}_{n-j}\leq \pi_A\,\omega_{1+\alpha}(t_n)\quad\text{for $1\le n\le N$.}
\end{equation}

Recall that the Mittag--Leffler function
$E_\alpha(z) = \sum_{k=0}^\infty\frac{z^k}{\Gamma(1+k\alpha)}$.
We state the following (slightly simplified) version of~\cite[Theorem~3.2]{LiaoMcLeanZhangZhao:2018}.
This result differs substantially from
the fractional Gr\"{o}nwall inequality of Jin et al. \cite[Theorem~4]{JinLiZhou:2017}
 since it is valid on very general nonuniform time meshes.
\begin{theorem}\label{thm: Gronwall}
Let \Ass{1}--\Ass{3} hold true. Suppose that the sequences $(\xi_1^n)_{n=1}^N$, $(\xi_2^n)_{n=1}^N$ are nonnegative.
Assume that $\lambda_0$ and $\lambda_1$ are non-negative constants and the maximum step size
$\tau\leq1/\sqrt[\alpha]{2\max\{1,\rho\}\pi_A\Gamma(2-\alpha)(\lambda_0+\lambda_1)}$.
If the nonnegative sequence~$(v^k)_{k=0}^N$ satisfies
\begin{equation}\label{eq: first Gronwall}
\sum_{k=1}^nA^{(n)}_{n-k}\diff v^k\le
	\lambda_{0}v^{n}+\lambda_1v^{n-1}+\xi_1^n+\xi_2^n\quad\text{for $1\le n\le N$,}
\end{equation}
then it holds that for $1\le n\le N$,
\begin{equation}\label{eq: Gronwall conclusion}
v^n\le2 E_\alpha\brab{2\max\{1,\rho\}\pi_A(\lambda_0+\lambda_1)t_n^\alpha}
	\biggl(v^0+\max_{1\le k\le n}\sum_{j=1}^k P^{(k)}_{k-j}\xi_1^j
+\pi_A\omega_{1+\alpha}(t_n)\max_{1\le j\le n}\xi_2^j\biggr).
\end{equation}
\end{theorem}

To facilitate our analysis, we now eliminate the historic term $H^{\ell}(t_{n})$ from
the fast L1  formula~\eqref{fastL1formulaNonuniform} for~$ \brat{\ffd{\alpha}u}^n$.
From the recurrence relationship \eqref{fastL1formulaRecurrence}, it is easy to see that
\begin{align*}
H^{\ell}(t_{k})=\sum_{j=1}^ke^{-\theta^{\ell}(t_{k}-t_{j})}b^{(j,\ell)}\nabla_{\tau}u^{j},\quad k\geq1, \; 1\leq \ell\leq N_q.
\end{align*}
Inserting this in \eqref{fastL1formulaNonuniform} and using
the definition \eqref{Coefficient1-fastL1formulaNonuniform}, one obtains the alternative formula
\begin{align}\label{fastL1formulaNonuniform-anotherform}
 \brat{\ffd{\alpha}u}^n
=&\,a_{0}^{(n)}\nabla_{\tau}u^n+\sum_{k=1}^{n-1}\frac{\nabla_{\tau}u^{k}}{\tau_k}
\int_{t_{k-1}}^{t_{k}}\sum_{\ell=1}^{N_q}\varpi^{\ell}e^{-\theta^{\ell}(t_{n}-s)}\zd s
= \sum_{k=1}^{n}A_{n-k}^{(n)}\nabla_{\tau}u^{k}, \quad n\geq1,
\end{align}
where the discrete convolution kernel $A_{n-k}^{(n)}$ is henceforth defined by
\begin{align}\label{Coefficient2-fastL1formulaNonuniform}
 A_{0}^{(n)}\defeq a_{0}^{(n)}, \quad
A_{n-k}^{(n)}\defeq\frac{1}{\tau_k}\int_{t_{k-1}}^{t_{k}}\sum_{\ell=1}^{N_q}\varpi^{\ell}e^{-\theta^{\ell}(t_{n}-s)}\zd s,
\quad 1\leq k\leq n-1,\; n\geq1.
\end{align}

The formula \eqref{fastL1formulaNonuniform-anotherform} takes the form of
\eqref{DiscreteCaputoDef}, and we now verify that our $A_{n-k}^{(n)}$
defined by \eqref{Coefficient2-fastL1formulaNonuniform} satisfy  \Ass{1} and \Ass{2},
allowing us to apply Theorem~\ref{thm: Gronwall} and establish the convergence of our computed solution.
Part (I) of the next lemma ensures that \Ass{1} is valid,
while part (II) implies that \Ass{2} holds true with $\pi_A=\frac{3}{2}$.

\begin{lemma}\label{lemma:Coefficient2-decreasingFastL1}
If the tolerance error $\epsilon$ of SOE satisfies
$\epsilon\leq\min\left\{\frac{1}{3}\omega_{1-\alpha}(T),\alpha\,\omega_{2-\alpha}(1)\right\}$,
then the discrete convolutional kernel $A_{n-k}^{(n)}$ of \eqref{Coefficient2-fastL1formulaNonuniform} satisfies
\begin{itemize}
\item [(I)] $ A_{k-1}^{(n)}>A_{k}^{(n)}>0,\;\;1\leq k\leq n-1;$\quad
(II) $ A_{0}^{(n)}=a_{0}^{(n)}$ and $A_{n-k}^{(n)}\geq \frac{2}{3}a_{n-k}^{(n)},\;\; 1\leq k\leq n-1.$
\end{itemize}
\end{lemma}
\begin{proof}
The definition \eqref{Coefficient1-L1formulaNonuniform} and Lemma \ref{lemma: L1 kernel} (i) yield
$$a_{0}^{(n)}-a_{1}^{(n)}>a_{0}^{(n)}-\omega_{1-\alpha}(\tau_n)
=\tfrac{\alpha}{\tau_n}\omega_{2-\alpha}(\tau_n)\geq \alpha\,\omega_{2-\alpha}(1)\geq\epsilon\,,$$
where the step size $\tau_n\leq 1$ and our hypothesis on $\epsilon$ are used.
The definition \eqref{Coefficient2-fastL1formulaNonuniform} and Lemma \ref{lemma:sum-of-exponentials} imply that $A_{0}^{(n)}=a_{0}^{(n)}>a_{1}^{(n)}+\epsilon>A_{1}^{(n)}$.
Lemma \ref{lemma:sum-of-exponentials} also shows that $\theta^{\ell}, \varpi^{\ell}>0$ for $\ell =1,\dots, N_q$;
the mean-value theorem now yields property (I). By Lemma~\ref{lemma: L1 kernel} (i)
and our hypothesis on $\epsilon$ we have
$\epsilon\leq\frac13\omega_{1-\alpha}(t_n - t_{k-1})<\frac1{3}a_{n-k}^{(n)}$ for $1\leq k\leq n-1$.
Hence Lemma \ref{lemma:sum-of-exponentials} gives $A_{n-k}^{(n)}\geq a_{n-k}^{(n)}-\epsilon\geq \frac2{3}a_{n-k}^{(n)}$ for
$1\leq k\leq n-1$. The proof is complete.
\end{proof}

\section{Global consistency error analysis}\label{sec:consistency}
\setcounter{equation}{0}

We now proceed with the consistency error analysis of our fast linearized method, and begin with the
consistency error of the standard L1 formula  $(\dfd{\alpha}u)^n$ of~\eqref{L1formulaNonuniform}.

\begin{lemma}\label{lemma:L1formulaNonuniform-consistence}
For $v\in C^2(0,T]$ with $\int_0^T t \,|v''(t)|\zd s < \infty$, one has
\begin{align*}
\absb{(\fd{\alpha}v)(t_n)-(\dfd{\alpha}v)^n}\leq
a_{0}^{(n)}G^n+\sum_{k=1}^{n-1}\brab{a_{n-k-1}^{(n)}-a_{n-k}^{(n)}}G^k,\quad n\geq1,
\end{align*}
where the L1 kernel $a_{n-k}^{(n)}$ is defined by \eqref{Coefficient1-L1formulaNonuniform} and $G^k :=2\int_{t_{k-1}}^{t_k}\bra{t-t_{k-1}}\abs{v''(t)}\zd t$.
\end{lemma}
\begin{proof}
From Taylor's formula with integral remainder, the truncation error of
the standard L1 formula at time $t=t_n$ is  (see \cite[Lemma 3.3]{LiaoLiZhangZhao:2017})
\begin{align}
(\fd{\alpha}v)(t_n)-(\dfd{\alpha}v)^n
=&\,\sum_{k=1}^n\int_{t_{k-1}}^{t_k}\omega_{1-\alpha}(t_n-s)\bra{v'(s)-\nabla_\tau v^k/\tau_k}\zd s\nonumber\\
=&\,\sum_{k=1}^n\int_{t_{k-1}}^{t_k}v''(t)\,\brab{\widetilde{\Pi_{1,k}}Q}(t)\zd{t}\,,\quad n\geq1,
\label{localApproximateError-integral}
\end{align}
where $Q(t)=\omega_{2-\alpha}(t_n-t)$ and we use the notation of the proof of Lemma~\ref{lemma: L1 kernel}.
By the error formula for linear interpolation \cite[Lemma 3.1]{LiaoLiZhangZhao:2017}, we have
\begin{align*}
\brab{\Err{1,k}Q}(t)=\int_{t_{k-1}}^{t_k}\chi_k(t,y)Q''(y)\zd{y},\quad t_{k-1}<t<t_k,\;1\leq k\leq n,
\end{align*}
where the Peano kernel $\chi_k(t,y)=\max\{t-y,0\}-\frac{t-t_{k-1}}{\tau_k}(t_{k}-y)$ satisfies
$$-\tfrac{t-t_{k-1}}{\tau_k}(t_{k}-t)\leq \chi_k(t,y)<0\qquad \text{for any $t,y\in(t_{k-1},t_k).$}$$
Observing that for each fixed $n\geq1$ the function $Q$ is decreasing and $Q''(t)=\omega_{-\alpha}(t_n-t)<0$,
we arrive at
the interpolation error $\brab{\Err{1,k}Q}(t)\geq0$ for $1\leq k\leq n$, with
\begin{align*}
\brab{\Err{1,n}Q}(t)
&\leq Q(t_{n-1})-\brab{\Pi_{1,n}Q}(t)=(t-t_{n-1})a_{0}^{(n)},\\
\brab{\Err{1,k}Q}(t) &\leq
(t_{k-1}-t)\int_{t_{k-1}}^{t_k}Q''(t)\zd{t}
\leq(t-t_{k-1})\kbrab{\omega_{1-\alpha}(t_n-t_{k})-\omega_{1-\alpha}(t_n-t_{k-1})}\nonumber\\
&\leq 2(t-t_{k-1})\brab{a_{n-k-1}^{(n)}-a_{n-k}^{(n)}},\qquad t\in(t_{k-1},t_{k}), \; 1\leq k\leq n-1,
\end{align*}
where Lemma \ref{lemma: L1 kernel} (ii) is used in the last inequality.
Thus, \eqref{localApproximateError-integral} yields
\begin{align*}
&\,\absb{(\fd{\alpha}v)(t_n)-(\dfd{\alpha}v)^n}
\leq\int_{t_{n-1}}^{t_n}\abs{v''(t)}\brab{\widetilde{\Pi_{1,n}}Q}(t)\zd{t}
+\sum_{k=1}^{n-1}\int_{t_{k-1}}^{t_k}\abs{v''(t)}\brab{\widetilde{\Pi_{1,k}}Q}(t)\zd{t}\\
&\,\hspace{1.4cm}\leq a_{0}^{(n)}\int_{t_{n-1}}^{t_n}\!\!(t-t_{n-1})\abs{v''(t)}\!\zd{t}
+2\sum_{k=1}^{n-1}\brab{a_{n-k-1}^{(n)}-a_{n-k}^{(n)}}\int_{t_{k-1}}^{t_k}\!\!(t-t_{k-1})\abs{v''(t)}\!\zd{t},
\end{align*}
and the desired result follows from the definition of $G^k$.
\end{proof}

\begin{remark}
Compared with the previous estimate in~\cite[Lemma~3.3]{LiaoLiZhangZhao:2017},
Lemma~\ref{lemma:L1formulaNonuniform-consistence} removes the time-step ratios restriction
$\rho_k\leq1$, which is an undesirable limitation on the mesh for the problems
that allow the rapid growth of  the solution at the time far away from $t=0$.
\end{remark}

We now focus on the fast L1 method by taking the initial singularity into account.
Here and hereafter, we denote $\hat{T}=\max\{1,T\}$ and $\hat{t}_{n}=\max\{1,t_{n}\}$ for $1\leq n\leq N$.

\begin{lemma}\label{lemma:fastL1formulaNonuniform-consistence}
Assume that $v\in C^2((0,T])$ and that there exists a constant $C_v>0$ such that
\begin{align}\label{regularityAssumption}
\absb{v'(t)}\leq C_v (1+t^{\sigma-1}),\quad \absb{v''(t)}\leq C_v (1+t^{\sigma-2}), \quad 0<t\leq T,
\end{align}
where $\sigma\in(0,1)\cup(1,2)$ is a parameter. Let $\Upsilon^j :=\brat{\fd{\alpha}v}(t_j)- \brat{\ffd{\alpha}v}^j$
denote the local consistency error of the fast L1 formula \eqref{fastL1formulaNonuniform-anotherform}.
Assume that the SOE tolerance error $\epsilon$ satisfies
$\epsilon\leq\frac{1}{3}\min\{\omega_{1-\alpha}(T),3\alpha\,\omega_{2-\alpha}(1)\}$.
Then the global consistency error
\begin{align}\label{globalfastL1ConsistenceError}
\sum_{j=1}^nP_{n-j}^{(n)}\absb{\Upsilon^j}
\leq&\,C_v\braB{\frac{\tau_1^{\sigma}}{\sigma}+
\frac{1}{1-\alpha}\max_{2\leq k\leq n}(t_{k}-t_1)^{\alpha}t_{k-1}^{\sigma-2}\tau_k^{2-\alpha}
+\frac{\epsilon}{\sigma}t_{n}^{\alpha}\hat{t}_{n-1}^{\,2}}
\end{align}
for $1\leq n\leq N$. Moreover, if the mesh satisfies \Ass{G}, then
\begin{align*}
\sum_{j=1}^nP_{n-j}^{(n)}\absb{\Upsilon^j}
\leq&\,\frac{C_v}{\sigma(1-\alpha)}\tau^{\min\{2-\alpha,\gamma\sigma\}}
+\frac{\epsilon}{\sigma}C_vt_{n}^{\alpha}\hat{t}_{n-1}^{\,2},\quad 1\leq n\leq N.
\end{align*}
\end{lemma}
\begin{proof}
The main difference between the fast L1 formula \eqref{fastL1formulaNonuniform-anotherform}
and the standard L1 formula \eqref{L1formulaNonuniform} is that the convolution kernel
is approximated by SOEs with an absolute  tolerance error $\epsilon$.
Thus, comparing the standard L1 formula \eqref{L1formulaNonuniform}
with the corresponding fast L1 formula  \eqref{fastL1formulaNonuniform-anotherform},
by  Lemma~\ref{lemma:sum-of-exponentials} and the regularity assumption \eqref{regularityAssumption} one has
\begin{align*}
\absb{ \brat{\ffd{\alpha}v}^j
-\brat{\dfd{\alpha}v}^j} &\leq \sum_{k=1}^{j-1}\frac{\absb{\nabla_{\tau}v^{k}}}{\tau_k}\int_{t_{k-1}}^{t_{k}}
\absB{\sum_{\ell=1}^{N_q}\varpi^{\ell}e^{-\theta^{\ell}(t_{j}-s)}
-\omega_{1-\alpha}(t_j-s)}\zd s,\\
&\leq \epsilon\sum_{k=1}^{j-1}\int_{t_{k-1}}^{t_{k}}\abs{v'(s)}\zd s
\leq C_v \brab{t_{j-1} + t_{j-1}^{\sigma}/\sigma}\epsilon\leq\frac{C_v}{\sigma}\hat{t}_{j-1}^{\,2}\epsilon, \quad j\geq1.
\end{align*}
Lemma \ref{lemma:sum-of-exponentials} implies that $\absb{A_{n-k}^{(n)}-a_{n-k}^{(n)}}\leq \epsilon$ for $1\leq k\leq n-1$.
Recalling that $A_{0}^{(n)}=a_{0}^{(n)}$, one has
$a_{j-k-1}^{(j)}-a_{j-k}^{(j)}\leq  A_{j-k-1}^{(j)}-A_{j-k}^{(j)}+2\epsilon$ for $1\leq k\leq j-1.$
Then Lemma \ref{lemma:L1formulaNonuniform-consistence} and the regularity assumption \eqref{regularityAssumption} lead to
\begin{align*}
\abs{(\fd{\alpha}v)(t_j)-\brat{\dfd{\alpha}v}^j}
\leq&\,
A_{0}^{(j)}G^j+\sum_{k=1}^{j-1}\brab{A_{j-k-1}^{(j)}-A_{j-k}^{(j)}}G^k+2\epsilon\sum_{k=1}^{j-1}G^k\\
\leq&\,
A_{0}^{(j)}G^j+\sum_{k=1}^{j-1}\brab{A_{j-k-1}^{(j)}-A_{j-k}^{(j)}}G^k
+4\epsilon\sum_{k=1}^{j-1}\int_{t_{k-1}}^{t_k}t\abs{v''(t)}\zd{t}\\
\leq&\,
A_{0}^{(j)}G^j+\sum_{k=1}^{j-1}\brab{A_{j-k-1}^{(j)}-A_{j-k}^{(j)}}G^k
+\frac{C_v}{\sigma}\hat{t}_{j-1}^{\,2}\epsilon,\quad j\geq1.
\end{align*}
Now a triangle inequality gives
\begin{align}\label{localfastL1ConsistenceError}
\absb{\Upsilon^j}
\leq&\,A_{0}^{(j)}G^j+\sum_{k=1}^{j-1}\brab{A_{j-k-1}^{(j)}-A_{j-k}^{(j)}}G^k
+ \frac{C_v}{\sigma}\hat{t}_{j-1}^{\,2}\epsilon,\quad j\geq1.
\end{align}
Multiplying the above inequality \eqref{localfastL1ConsistenceError} by
$P_{n-j}^{(n)}$ and summing the index $j$ from $1$ to $n$, one can
exchange the order of summation and apply
the definition \eqref{discreteConvolutionKernel-RL} of $P_{n-j}^{(n)}$ to obtain
\begin{align}\label{globalApproximateError-immediate}
\sum_{j=1}^nP_{n-j}^{(n)}&\,\absb{\Upsilon^j}\leq
\sum_{j=1}^nP_{n-j}^{(n)}A_0^{(j)}G^j+
\sum_{j=2}^nP_{n-j}^{(n)}\sum_{k=1}^{j-1}\brab{A_{j-k-1}^{(j)}-A_{j-k}^{(j)}}G^k
+C_v\frac{\epsilon}{\sigma}\sum_{j=2}^nP_{n-j}^{(n)}\hat{t}_{j-1}^{\,2}\nonumber\\
=&\,\sum_{j=1}^nG^jP_{n-j}^{(n)}A_0^{(j)}
+\sum_{k=1}^{n-1}G^k\sum_{j=k+1}^nP_{n-j}^{(n)}\brab{A_{j-k-1}^{(j)}-A_{j-k}^{(j)}}
+C_v\hat{t}_{n-1}^{\,2}\frac{\epsilon}{\sigma}\sum_{j=2}^nP_{n-j}^{(n)}\nonumber\\
\leq&\,  \sum_{k=1}^{n}P_{n-k}^{(n)}A_{0}^{(k)}G^k+\sum_{k=1}^{n-1}P_{n-k}^{(n)}A_{0}^{(k)}G^k
+\frac{C_v}{\sigma}t_{n}^{\alpha}\hat{t}_{n-1}^{\,2}\epsilon,
\end{align}
where the property \eqref{eq: P bound} with $\pi_A=3/2$ is used in the last inequality.
If the SOE approximation error $\epsilon\leq\frac{1}{3}\min\{\omega_{1-\alpha}(T),3\alpha\,\omega_{2-\alpha}(1)\}\,,$
Lemma \ref{lemma:Coefficient2-decreasingFastL1} (II) and
Lemma \ref{lemma: L1 kernel} (i) imply that
$A_{0}^{(k)}=a_{0}^{(k)}=\omega_{2-\alpha}(\tau_k)/\tau_k$,
$A_{k-2}^{(k)}\geq \tfrac2{3}a_{k-2}^{(k)}\geq\tfrac2{3}\omega_{1-\alpha}(t_k-t_1)$, and then
\begin{align*}
A_{0}^{(k)}/{A_{k-2}^{(k)}}\leq\tfrac{3}{2(1-\alpha)}(t_k-t_1)^{\alpha}\tau_k^{-\alpha},\quad 2\leq k\leq n\le N.
\end{align*}
Furthermore, the identical property~\eqref{eq: P A} for the complementary kernel~$P^{(n)}_{n-j}$ gives
$$P_{n-1}^{(n)}A_{0}^{(1)}\leq1\quad \mbox{and}\quad
\sum_{k=2}^{n-1}P_{n-k}^{(n)}A_{k-2}^{(k)}\leq\sum_{k=2}^{n}P_{n-k}^{(n)}A_{k-2}^{(k)}=1.$$
The regularity assumption \eqref{regularityAssumption} gives
$G^1\leq C_v\tau_1^{\sigma}/\sigma$ and $G^k\leq C_vt_{k-1}^{\sigma-2}\tau_k^2$ $(2\leq k\leq n).$
Thus it follows from \eqref{globalApproximateError-immediate}  that
\begin{align*}
\sum_{j=1}^nP_{n-j}^{(n)}\absb{\Upsilon^j}\leq&\,2G^1+2\sum_{k=2}^{n}P_{n-k}^{(n)}A_{0}^{(k)}G^k
+\frac{C_v}{\sigma}t_{n}^{\alpha}\hat{t}_{n-1}^{\,2}\epsilon\\
\leq&\,
C_v\frac{\tau_1^{\sigma}}{\sigma}+\frac{C_v}{1-\alpha}\sum_{k=2}^{n}P_{n-k}^{(n)}A_{k-2}^{(k)}
(t_{k}-t_1)^{\alpha}t_{k-1}^{\sigma-2}\tau_k^{2-\alpha}
+\frac{C_v}{\sigma}t_{n}^{\alpha}\hat{t}_{n-1}^{\,2}\epsilon\\
\leq&\,
C_v\braB{\frac{\tau_1^{\sigma}}{\sigma}
+\frac{1}{1-\alpha}\max_{2\leq k\leq n}(t_{k}-t_1)^{\alpha}t_{k-1}^{\sigma-2}\tau_k^{2-\alpha}
+\frac{1}{\sigma}t_{n}^{\alpha}\hat{t}_{n-1}^{\,2}\epsilon},\quad 1\leq n\leq N.
\end{align*}
The claimed estimate \eqref{globalfastL1ConsistenceError} is verified.
In particular, if \Ass{G} holds, one has
\begin{align*}
t_{k}^{\alpha}t_{k-1}^{\sigma-2}\tau_k^{2-\alpha}
\leq&\, C_{\gamma}t_k^{\sigma-2+\alpha}\tau_k^{2-\alpha-\beta}\tau^{\beta}\min\{1,t_k^{\beta-\beta/\gamma}\}\\
\leq&\, C_{\gamma}t_k^{\sigma-\beta/\gamma}\brab{\tau_k/t_k}^{2-\alpha-\beta}\tau^{\beta}
\leq C_{\gamma}t_k^{\max\{0,\sigma-(2-\alpha)\gamma\}}\tau^{\beta},\quad 2\leq k\leq N,
\end{align*}
where $\beta=\min\{2-\alpha,\gamma\sigma\}$. The final estimate follows since
$\tau_1^{\sigma}\leq C_{\gamma}\tau^{\gamma\sigma}\leq C_{\gamma}\tau^{\beta}$.
\end{proof}


Next lemma describes the global consistency error of Newton's linearized approach, which is
smaller than that generated by the above L1 approximation. In addition,
there is no error in the linearized approximation if $f=f(u)$ is a linear function.
\begin{lemma}\label{lemma:NewtonlinearizedConsistence}
Assume that $v\in C([0,T])\cap C^2((0,T])$ satisfies the regularity condition \eqref{regularityAssumption},
and the nonlinear function $f=f(u)\in C^2(\mathbb{R})$.
Denote $v^n=v(t_n)$ and the local truncation error
$\mathcal{R}_f^n=f(v^n)-f(v^{n-1})-f'(v^{n-1})\nabla_{\tau}v^{n}$ such that the global consistency error
\begin{align*}
\sum_{j=1}^nP_{n-j}^{(n)}\absb{\mathcal{R}_f^j}\leq C_v\tau_1^{\alpha}\bra{\tau_1^{2}+\tau_1^{2\sigma}/\sigma^2}
+C_vt_n^{\alpha}\max_{2\leq j\leq n} \brab{\tau_j^2+t_{j-1}^{2\sigma-2}\tau_j^{2}},\quad 1\leq n\leq N.
\end{align*}
Moreover, if the assumption \Ass{G} holds, one has
\begin{align*}
\sum_{j=1}^nP_{n-j}^{(n)}\absb{\mathcal{R}_f^j}
\leq&\,C_v\tau^{\min\{2,2\gamma\sigma\}}\max\{1,\tau^{\gamma\alpha}/\sigma^2\},\quad 1\leq n\leq N.
\end{align*}
\end{lemma}
\begin{proof}
Applying the formula of
Taylor expansion with integral remainder, one has
\begin{align*}
\mathcal{R}_f^j=(\nabla_{\tau}v^{j})^2\int_0^1f''\brab{v^{j-1}+s\nabla_{\tau}v^{j}}(1-s)\zd{s},\quad j\geq1.
\end{align*}
Under the regularity conditions, one has
$\absb{\mathcal{R}_f^{1}}\leq C_v\brab{\int_{t_{0}}^{t_1}\abs{v'(t)}\!\zd{t}}^2
\leq C_v\bra{\tau_1^{2}+\tau_1^{2\sigma}/\sigma^2},$
\begin{align*}
\absb{\mathcal{R}_f^{j}}\leq C_v\braB{\int_{t_{j-1}}^{t_j}\abs{v'(t)}\!\zd{t}}^2
\leq C_v\brab{\tau_j^2+t_{j-1}^{2\sigma-2}\tau_j^{2}},\quad 2\leq j\leq N.
\end{align*}
Note that, Lemma \ref{lemma:Coefficient2-decreasingFastL1} (II) and
the definition \eqref{Coefficient1-L1formulaNonuniform} give
$A_{0}^{(k)}=a_{0}^{(k)}=\omega_{2-\alpha}(\tau_k)/\tau_k$,
so the identical property \eqref{eq: P A} shows
$P_{n-1}^{(n)}\leq1/A_{0}^{(1)}\leq\Gamma(2-\alpha)\tau_1^{\alpha}$.
Moreover, the bounded estimate \eqref{eq: P bound} with $\pi_A=\frac32$
gives $\sum_{j=2}^nP_{n-j}^{(n)}\leq \frac32\omega_{1+\alpha}(t_n)$.
Thus, it follows that
\begin{align*}
\sum_{j=1}^nP_{n-j}^{(n)}\absb{\mathcal{R}_f^j}\leq&\,P_{n-1}^{(n)}\absb{\mathcal{R}_f^1}
+\sum_{j=2}^nP_{n-j}^{(n)}\absb{\mathcal{R}_f^j}
\leq
C_v\tau_1^{\alpha}\absb{\mathcal{R}_f^1}+C_vt_n^{\alpha}\max_{2\leq j\leq n}\absb{\mathcal{R}_f^j}\\
\leq&\,
C_v\tau_1^{\alpha}\bra{\tau_1^{2}+\tau_1^{2\sigma}/\sigma^2}
+C_vt_n^{\alpha}\max_{2\leq j\leq n} \brab{\tau_j^2+t_{j-1}^{2\sigma-2}\tau_j^{2}},\quad 1\leq n\leq N.
\end{align*}
If \Ass{G} holds, one has
$\tau_j^2\leq C_{\gamma}\tau^{2}\min\{1,t_j^{2-2/\gamma}\}\leq C_{\gamma}\tau^{\beta}\min\{1,t_j^{2-2/\gamma}\}$, and
\begin{align*}
t_{j-1}^{2\sigma-2}\tau_j^{2}
\leq&\, C_{\gamma}t_j^{2\sigma-2}\tau_j^{2-\beta}\tau^{\beta}\min\{1,t_j^{\beta-\beta/\gamma}\}\\
\leq&\, C_{\gamma}t_j^{2\sigma-\min\{2,2\gamma\sigma\}/\gamma}\brab{\tau_k/t_k}^{2-\beta}\tau^{\beta}
\leq C_{\gamma}t_k^{\max\{0,2\sigma-2/\gamma\}}\tau^{\beta},\quad 2\leq j\leq N,
\end{align*}
where $\beta=\min\{2,2\gamma\sigma\}$.
The second estimate follows since
$\tau_1^{2\sigma}\leq C_{\gamma}\tau^{2\gamma\sigma}\leq C_{\gamma}\tau^{\beta}$.
\end{proof}

\section{Unconditional convergence}\label{sec:error}
\setcounter{equation}{0}

Assume that the time mesh fulfills \Ass{3} and \Ass{G} in the error analysis.
We improve the discrete $H^2$ energy method in
\cite{LiaoSun:2010,LiaoSunShi:2010,LiaoSunShi:2010B}
to prove the unconditional convergence of discrete
solution to the two-level linearized scheme \eqref{diffScheme}.
In this section, $K_0$, $\tau_0$, $\tau_1$, $\tau_0^{*}$, $h_0$, $\epsilon_0$ and
any numeric subscripted $c$, such as $c_0$, $c_1$, $c_2$ and so on,
are fixed values, which are always dependent on the given data and the solution,
but independent of the time-space grid steps and the inductive index $k$ in the mathematical induction as well.
To make our ideas more clearly, four steps to obtain unconditional error estimate are listed in four subsections.

\subsection{STEP 1: construction of coupled discrete system}
We introduce a function
$w\defeq\fd{\alpha}u-f(u)$ with the initial-boundary values
$w(\vecx ,0)\defeq\Delta u^0(\vecx )$ for $\vecx \in \Omega$ and $w(\vecx ,t)\defeq-f(0)$ for $\vecx \in \partial \Omega$.
The problem \eqref{Problem-1} can be formulated into
\begin{align*}
w=&\,\fd{\alpha}u-f(u),
\quad \vecx \in \bar{\Omega},\; 0<t\leq T;\\
w=&\,\Delta u,\quad \vecx \in \Omega,\; 0\leq t\leq T.
\end{align*}
Let $w_{h}^n$ be the numerical approximation of function $W_{h}^n=w(\vecx_h,t_n)$
for $\vecx_h\in\bar{\Omega}_h$. As done in subsection 2.2,
one has an auxiliary discrete system: to seek $\{u_h^{n},\,w_h^n\}$ such that
\begin{align}
w_h^n=&\,\brat{\ffd{\alpha}u_h}^n-f(u_h^{n-1})-f'(u_h^{n-1})\nabla_{\tau}u_h^{n}\,,
\quad  \vecx_h\in\bar{\Omega}_h, \; 1\leq n\leq N;\label{differenceScheme-auxiliary1}\\
w_h^n=&\,\Delta_hu_h^n\,,
\quad  \vecx_h\in\Omega_h, \; 0\leq n\leq N;\label{differenceScheme-auxiliary2}\\
u_h^0=&\,u^0(\vecx_h),\;\; \vecx_h\in\bar{\Omega}_h\,;\quad
u_h^n=0,\; \;\vecx_h\in\partial{\Omega}_h\,,1\leq n\leq N.\label{differenceScheme-auxiliary3}
\end{align}
Obviously, by eliminating the auxiliary function $w_h^n$ in above discrete system,
one directly arrives at the computational scheme \eqref{diffScheme}.
Alternately, the solution properties of two-level linearized method
\eqref{diffScheme} can be studied via the auxiliary discrete system
\eqref{differenceScheme-auxiliary1}-\eqref{differenceScheme-auxiliary3}.

\subsection{STEP 2: reduction of coupled error system}
Let $\tilde{u}_{h}^n=U_{h}^n-u_{h}^n$, $\tilde{w}_{h}^n=W_{h}^n-w_{h}^n$ be the solution errors for $\vecx_h\in\bar{\Omega}_h$.
We now have an error system with respect to the error function $\{\tilde{w}_h^n\}$ as
\begin{align}
\tilde{w}_h^n=&\,\brat{\ffd{\alpha}\tilde{u}_h}^n-\mathcal{N}_h^n+\xi_h^n\,,
\quad  \vecx_h\in\bar{\Omega}_h, \; 1\leq n\leq N;\label{errorSystem-auxiliary1}\\
\tilde{w}_h^n=&\,\Delta_h\tilde{u}_h^n+\eta_h^n\,,
\quad  \vecx_h\in\Omega_h, \; 0\leq n\leq N;\label{errorSystem-auxiliary2}\\
\tilde{u}_h^0=&\,0,\;\; \vecx_h\in\bar{\Omega}_h\,;\quad
\tilde{u}_h^n=0,\;\; \vecx_h\in\partial{\Omega}_h\,,1\leq n\leq N,\label{errorSystem-auxiliary3}
\end{align}
where $\xi_h^n$ and $\eta_h^n$ denote  temporal and spatial truncation errors, respectively, and
\begin{align}\label{nonlinearTerm-errorSystem}
\mathcal{N}_h^n\defeq&\,f'(u_h^{n-1})\nabla_{\tau}\tilde{u}_h^{n}+f(U_h^{n-1})-f(u_h^{n-1})
+\bra{f'(U_h^{n-1})-f'(u_h^{n-1})}\nabla_{\tau}U_h^{n}\nonumber\\
=&\,f'(u_h^{n-1})\nabla_{\tau}\tilde{u}_h^{n}
+\tilde{u}_h^{n-1}\int_0^1f'\brab{s U_h^{n-1}+(1-s)u_h^{n-1}}\zd{s}\nonumber\\
&\,\hspace{2.4cm}+\tilde{u}_h^{n-1}\nabla_{\tau}U_h^{n}\int_0^1f''\brab{s U_h^{n-1}+(1-s)u_h^{n-1}}\zd{s}\,.
\end{align}
Acting the difference operators $\Delta_h$ and $\ffd{\alpha}$ on the equations
\eqref{errorSystem-auxiliary1}-\eqref{errorSystem-auxiliary2}, respectively, gives
 \begin{align*}
\Delta_h\tilde{w}_h^n=&\,\brat{\ffd{\alpha}\Delta_h\tilde{u}_h}^n-\Delta_h\mathcal{N}_h^n+\Delta_h\xi_h^n\,,
\quad  \vecx_h\in\Omega_h, \; 1\leq n\leq N;\\
\brat{\ffd{\alpha}\tilde{w}_h}^n=&\,\brat{\ffd{\alpha}\Delta_h\tilde{u}_h}^n+\brat{\ffd{\alpha}\eta_h}^n\,,
\quad  \vecx_h\in\Omega_h, \; 1\leq n\leq N.
\end{align*}
By eliminating the term $\brat{\ffd{\alpha}\Delta_h\tilde{u}_h}^n$ in the above two equations, one gets
 \begin{align}
\brat{\ffd{\alpha}\tilde{w}_h}^n=&\,
\Delta_h\tilde{w}_h^n+\Delta_h\mathcal{N}_h^n+\brat{\ffd{\alpha}\eta_h}^n-\Delta_h\xi_h^n
\quad  \vecx_h\in\Omega_h, \; 1\leq n\leq N;\label{errorSystem-auxiliary1OPERATED}\\
\tilde{w}_h^0=&\,\eta_h^0,\;\; \vecx_h\in\bar{\Omega}_h\,;\quad
\tilde{w}_h^n=0,\;\; \vecx_h\in\partial{\Omega}_h\,,1\leq n\leq N;\label{errorSystem-auxiliary2OPERATED}
\end{align}
where the initial and boundary conditions are derived from the error system
\eqref{errorSystem-auxiliary1}-\eqref{errorSystem-auxiliary3}.

\subsection{STEP 3: continuous analysis of truncation error}
According to the first regularity condition in \eqref{SolutionRegularityAssumption}, one has
\begin{align}
\mynormb{\eta^n}\leq
c_1h^2,\quad 0\leq n\leq N.\label{globalSpaceError-1}
\end{align}
Since the spatial error $\eta_h^n$ is defined uniformly
at the time $t=t_n$ (there is no temporal error in the equation \eqref{differenceScheme-auxiliary2}),
we can define a continuous function $\eta_{h}(t)$ for $\vecx_h=(x_i,y_j)\in\Omega_h,$
\begin{align*}
\eta_{h}(t)=&\,\frac{h_1^2}6\int_0^1\kbrab{\partial_{x}^{(4)}u(x_i-s h_1,y_j,t)
+\partial_{x}^{(4)}u(x_i+s h_1,y_j,t)}(1-s)^3\zd{s}\\
&\,+\frac{h_2^2}6\int_0^1\kbrab{\partial_{y}^{(4)}u(x_i,y_j-s h_2,t)
+\partial_{y}^{(4)}u(x_i,y_j+s h_2,t)}(1-s)^3\zd{s}\,,
\end{align*}
such that $\eta_h^n=\eta_h(t_n)$.
The second condition in \eqref{SolutionRegularityAssumption} implies $\mynorm{\eta'(t)}\leq C_uh^2(1+t^{\sigma-1})$.
Hence, applying the fast L1 formula \eqref{fastL1formulaNonuniform-anotherform} and
the equality \eqref{eq: P A}, one has
\begin{align}
\sum_{j=1}^nP_{n-j}^{(n)}\mynormb{\brat{\ffd{\alpha}\eta}^j}\leq&\,
\sum_{j=1}^nP_{n-j}^{(n)}\sum_{k=1}^jA_{j-k}^{(j)}\mynormb{\nabla_{\tau}\eta^k}
=\sum_{k=1}^n\mynormb{\nabla_{\tau}\eta^k}\leq \frac{c_2}{\sigma}\hat{t}_n^{\,2}h^2.
\label{globalSpaceError-2}
\end{align}

Since the time truncation error $\xi_h^n$ in \eqref{errorSystem-auxiliary1} is defined
uniformly with respect to grid point $\vecx_h\in\bar\Omega_h$,
we can define a continuous function $\xi^n(\vecx )=\xi_1^n(\vecx )+\xi_2^n(\vecx )$, where
$\xi_1^n$, $\xi_2^n$ denotes the truncation errors of
fast L1 formula and Newton's linearized approach respectively,
\begin{align*}
\xi_1^n=\brat{\fd{\alpha}u}(t_n)-\brat{\ffd{\alpha}u}^n,\quad
\xi_2^n=\brab{\nabla_{\tau}u(t_n)}^2\int_0^1f''\brab{u(t_{n-1})
+s\nabla_{\tau}u(t_{n})}(1-s)\zd{s},
\end{align*}
such that $\xi_h^n=\xi^n(x_i,y_j)$ for $\vecx_h\in\bar\Omega_h$.
By the Taylor expansion formula, one has
\begin{align*}
\Delta_h\brab{\xi_{1}^n}_{ij}=&\,\int_0^1\kbrab{\partial_{xx}\xi_1^n(x_i-s h_1,y_j)
+\partial_{xx}\xi_1^n(x_i+s h_1,y_j)}(1-s)\zd{s}\\
&\,+\int_0^1\kbrab{\partial_{yy}\xi_1^n(x_i,y_j-s h_2)
+\partial_{yy}\xi_1^n(x_i,y_j+s h_2)}(1-s)\zd{s}\,,\quad 1\leq n\leq N.
\end{align*}
Applying Lemma \ref{lemma:fastL1formulaNonuniform-consistence} with the second
and third regularity conditions in \eqref{SolutionRegularityAssumption}, we have
\begin{align*}
\sum_{j=1}^nP_{n-j}^{(n)}\mynormb{\Delta_h\xi_1^j}
\leq&\,\frac{C_u}{\sigma(1-\alpha)}\tau^{\min\{2-\alpha,\gamma\sigma\}}
+\frac{C_u}{\sigma}t_{n}^{\alpha}\hat{t}_{n-1}^{\,2}\epsilon,\quad 1\leq n\leq N.
\end{align*}
Similarly, one can write out an integral expression of $\Delta_h\brab{\xi_{2}^n}_{ij}$
by using the Taylor expansion.
Assuming $f\in C^4(\mathbb{R})$ and
taking $\tau\leq \tau_1=\sqrt[\gamma\alpha]{\sigma}$ such that
$\tau^{\gamma\alpha}\leq \tau_1^{\gamma\alpha}=\sigma$,
we apply Lemma \ref{lemma:NewtonlinearizedConsistence} with the second regularity
condition in \eqref{SolutionRegularityAssumption} to find,
\begin{align*}
\sum_{j=1}^nP_{n-j}^{(n)}\mynormb{\Delta_h\xi_2^j}
\leq C_u\tau^{\min\{2,2\gamma\sigma\}}\max\{1,\tau^{\gamma\alpha}/\sigma^2\}
\leq \frac{C_u}{\sigma}\tau^{\min\{2,2\gamma\sigma\}},\quad 1\leq n\leq N.
\end{align*}
Thus, the triangle inequality leads to
\begin{align}
\sum_{j=1}^nP_{n-j}^{(n)}\mynormb{\Delta_h\xi^j}\leq&\,
\frac{c_3}{\sigma(1-\alpha)}\tau^{\min\{2-\alpha,\gamma\sigma\}}
+\frac{c_4}{\sigma}t_{n}^{\alpha}\hat{t}_{n-1}^{\,2}\epsilon\,,
\quad 1\leq n\leq N.\label{globalTimeError-gradedGrid}
\end{align}

\subsection{STEP 4: error estimate by mathematical induction}

For a positive constant $C_0$, let $\mathcal{B}(0,C_0)$ be a ball in the space
of grid functions on $\bar\Omega_h$ such that
$\max\big\{\mynormt{\psi}_{\infty},\mynormt{\nabla_h\psi},\mynormt{\Delta_h\psi}\big\}\leq C_0$
for any grid function $\{\psi_h\}\in\mathcal{B}(0,C_0)$.
Always, we need the following result to treat the nonlinear terms but leave the proof to Appendix A.
\begin{lemma}\label{lemma:BoundednessOf2rdDiscretederivate}
 Let $F\in C^2(\mathbb{R})$ and a grid function $\{\psi_h\}\in\mathcal{B}(0,C_0)$.
 Thus there is a constant $C_F>0$ dependent on $C_0$ and $C_{\Omega}$ such that,
 $\mynorm{\Delta_h\kbra{F(\psi)v}}\leq C_F\mynorm{\Delta_hv}$ for any $\{v_h\}\in\mathcal{V}_{h}$.
\end{lemma}

Under the regularity assumption \eqref{SolutionRegularityAssumption} with $U_h^k=u(\vecx_h,t_k)$, we define a constant
$$K_0=\frac13\max_{0\leq k\leq N}\big\{\mynormb{U^k}_{\infty},\mynormb{\nabla_hU^k},\mynormb{\Delta_hU^k}\big\}.$$
For a smooth function $F\in C^2(\mathbb{R})$ and any grid function $\{v_h\}\in\mathcal{V}_{h}$,
we denote the maximum value of $C_F$ in Lemma \ref{lemma:BoundednessOf2rdDiscretederivate} as $c_0$ such that
\begin{align}\label{BoundednessOf2rdDerivate}
\mynorm{\Delta_h\kbra{F(w)v}}\leq c_0\mynorm{\Delta_hv}\quad
\text{for any grid function $\{w_h\}\in\mathcal{B}(0,K_0+1)$.}
\end{align}

Let $c_5$ be the maximum value of $C_{\Omega}$ to verify the embedding inequalities
in \eqref{discreteNormsEmbedding}, and
\begin{align*}
c_6=\max\{1,c_5\}E_{\alpha}\brab{3\max\{1,\rho\}(2K_0+3)c_0T^{\alpha}},\quad
c_7=3c_1+\frac{2c_2}{\sigma}\hat{T}^{2}+3(2K_0+3)c_0c_1T^{\alpha}.
\end{align*}
Also let $\tau_0^{*}=1/\sqrt[\alpha]{3\max\{1,\rho\}\Gamma(2-\alpha)(2K_0+3)c_0}\,,$ and
\begin{align*}
\tau_0=\sqrt[\gamma\alpha]{\frac{\sigma(1-\alpha)}{6c_3c_6}}\,,\quad
h_0=\frac{1}{\sqrt{3c_6c_7}}\,,\quad
\epsilon_0=\min\Big\{\frac{\sigma }{6c_4c_6\hat{T}^{2}T^{\alpha}},\frac{1}{3}\omega_{1-\alpha}(T),\alpha\,\omega_{2-\alpha}(1)\Big\}.
\end{align*}
For the simplicity of presentation, denote 
\begin{align*}
&E_k\defeq E_{\alpha}\brab{3\max\{1,\rho\}(2K_0+3)c_0t_k^{\alpha}},\\
&\mathcal{T}^k\defeq\frac{2c_3}{\sigma(1-\alpha)}\tau^{\min\{2-\alpha,\gamma\sigma\}}
+\braB{2c_1+\frac{2c_2}{\sigma}\hat{t}_k^{\,2}
+3(2K_0+3)c_0c_1t_k^{\alpha}}h^2+\frac{2c_4}{\sigma}t_{k}^{\alpha}\hat{t}_{k-1}^{\,2}\epsilon\,,
\end{align*}
where $1\leq k\leq N$. We now apply the mathematical induction to prove that
\begin{align}\label{ErrorEstimate}
\mynormb{\Delta_h\tilde{u}^k}\leq E_k\mathcal{T}^k+c_1h^2
\quad\text{for $1\leq k\leq N$},
\end{align}
if the time-space grids and the SOE approximation satisfies
\begin{align}\label{ErrorEstimate-condition}
\tau\leq \min\{\tau_0,\tau_1,\tau_0^{*}\}, \quad h\leq h_0,\quad \epsilon\leq\epsilon_0.
\end{align}
Note that, the restrictions in \eqref{ErrorEstimate-condition} ensures
the error function $\{\tilde{u}_h^{k}\}\in\mathcal{B}(0,1)$ for $1\leq k\leq N$.


Consider $k=1$ firstly. Since $\tilde{u}_h^0=0$, $\{u_h^0\}\in\mathcal{B}(0,K_0)\subset\mathcal{B}(0,K_0+1)$ and
the nonlinear term \eqref{nonlinearTerm-errorSystem} gives
$\mathcal{N}_h^1=f'(u_h^{0})\tilde{u}_h^{1}$.
For the function $f\in C^3(\mathbb{R})$, the inequality \eqref{BoundednessOf2rdDerivate} implies
\begin{align}\label{nonlinearTerm-estimateFirstLevel}
\mynormb{\Delta_h\mathcal{N}^1}=\mynormt{\Delta_h\bra{f'(u^{0})\tilde{u}^{1}}}
\leq c_0\mynormb{\Delta_h\tilde{u}^{1}}\leq c_0\mynormb{\tilde{w}^{1}}+c_0c_1h^2\,,
\end{align}
where the equation \eqref{errorSystem-auxiliary2} and the estimate \eqref{globalSpaceError-1} are used.
Taking the inner product of the equation \eqref{errorSystem-auxiliary1OPERATED} (for $n=1$) by $\tilde{w}_h^{1}$, one
gets
\begin{align*}
A_0^{(1)}\myinnerb{\nabla_{\tau}\tilde{w}^{1},\tilde{w}^{1}}
\leq\myinnerb{\Delta_h\mathcal{N}^1,\tilde{w}^{1}}
+\myinnerb{\brat{\ffd{\alpha}\eta}^1-\Delta_h\xi^1,\tilde{w}^{1}},
\end{align*}
because the zero-valued boundary condition in \eqref{errorSystem-auxiliary2OPERATED}
leads to $\myinnerb{\Delta_h\tilde{w}^{1},\tilde{w}^{1}}\leq0$.
With the view of Cauchy-Schwarz inequality and \eqref{nonlinearTerm-estimateFirstLevel}, one has
$\myinnerb{\nabla_{\tau}\tilde{w}^{1},\tilde{w}^{1}}\geq \mynormb{\tilde{w}^{1}}\nabla_{\tau}\brab{\mynormb{\tilde{w}^{1}}}$
and then
\begin{align*}
A_0^{(1)}\nabla_{\tau}\brab{\mynormb{\tilde{w}^{1}}}
\leq&\mynormb{\Delta_h\mathcal{N}^1}
+\mynormb{\brat{\ffd{\alpha}\eta}^1-\Delta_h\xi^1}
\leq c_0\mynormb{\tilde{w}^{1}}
+\mynormb{\brat{\ffd{\alpha}\eta}^1-\Delta_h\xi^1}+c_0c_1h^2\,.
\end{align*}
Setting $\tau_1\leq\tau_0^{*}\leq 1/\sqrt[\alpha]{3\max\{1,\rhomax\}\Gamma(2-\alpha)c_0}$, we apply Theorem \ref{thm: Gronwall}
(discrete fractional Gr\"{o}nwall inequality) with $\xi_1^1=\mynormb{\brat{\ffd{\alpha}\eta}^1-\Delta_h\xi^1}$ and $\xi_2^1=c_0c_1h^2$ to get
\begin{align*}
\mynormb{\tilde{w}^{1}}\leq&  E_\alpha\brab{3\max\{1,\rhomax\}c_0t_1^\alpha}
\braB{2\mynormb{\eta^{0}}
+2P_0^{(1)}\mynormb{\brat{\ffd{\alpha}\eta}^1-\Delta_h\xi^1}+3c_0c_1\omega_{1+\alpha}(t_1)h^2}\\
\leq&E_1\braB{\frac{2c_3}{\sigma(1-\alpha)}\tau^{\min\{2-\alpha,\gamma\sigma\}}
+2c_1h^2+\frac{2c_2}{\sigma}\hat{t}_1^{\,2}h^2+3c_0c_1\omega_{1+\alpha}(t_1)h^2}\leq
E_1\mathcal{T}^1,
\end{align*}
where the initial condition \eqref{errorSystem-auxiliary2OPERATED} and the error estimates \eqref{globalSpaceError-1}-\eqref{globalTimeError-gradedGrid} are used.
Thus, the equation \eqref{errorSystem-auxiliary2} and the inequality \eqref{globalSpaceError-1} yield
the estimate \eqref{ErrorEstimate} for $k=1$,
\begin{align*}
\mynormb{\Delta_h\tilde{u}^1}\leq\mynormb{\tilde{w}^{1}}+\mynormb{\eta^{1}}
\leq E_1\mathcal{T}^1+c_1h^2\,.
\end{align*}

Assume that the error estimate \eqref{ErrorEstimate} holds for $1\leq k\leq n-1$ ($n\geq2$).
Thus we apply the embedding inequalities in \eqref{discreteNormsEmbedding} to get
\begin{align*}
\max\big\{\mynormb{\tilde{u}^{k}}_{\infty},\mynormb{\nabla_h\tilde{u}^{k}},\mynormb{\Delta_h\tilde{u}^{k}}\big\}\leq
\max\{1,c_5\}\brab{E_k\mathcal{T}^k+c_1h^2},\quad1\leq k\leq n-1.
\end{align*}
Under the priori settings in \eqref{ErrorEstimate-condition}, we have the error function
$\{\tilde{u}_h^{k}\}\in\mathcal{B}(0,1)$, the discrete solution $\{u_h^{k}\}\in\mathcal{B}(0,K_0+1)$
for $1\leq k\leq n-1$, and the continuous solution $\{U_h^{k}\}\in\mathcal{B}(0,K_0)\subset\mathcal{B}(0,K_0+1)$.
Then, for the function $f\in C^4(\mathbb{R})$,
one applies the inequality \eqref{BoundednessOf2rdDerivate} to find that
\begin{align*}
&\mynormb{\Delta_h\!\kbra{f'(u^{n-1})\nabla_{\tau}\tilde{u}^{n}}}
\leq c_0\mynormb{\Delta_h\nabla_{\tau}\tilde{u}^{n}}
\leq c_0\mynormb{\Delta_h\tilde{u}^{n}}+c_0\mynormb{\Delta_h\tilde{u}^{n-1}}\,,\\
&\mynormb{\Delta_h\!\kbra{\tilde{u}^{n-1}f'\brab{s U^{n-1}+(1-s)u^{n-1}}}}
\leq c_0\mynormb{\Delta_h\tilde{u}^{n-1}}\,,\\
&\mynormb{\Delta_h\!\kbra{\tilde{u}^{n-1}\nabla_{\tau}U^{n}f''\brab{s U^{n-1}+(1-s)u^{n-1}}}}
\leq c_0 \mynormb{\Delta_h\brat{\tilde{u}^{n-1}\nabla_{\tau}U^{n}}}
\leq 2c_0K_0\mynormb{\Delta_h\tilde{u}^{n-1}},
\end{align*}
where $0\leq s\leq1$. From the expression \eqref{nonlinearTerm-errorSystem}
of $\mathcal{N}^n$ and the triangle inequality, one has
\begin{align}\label{nonlinearTerm-estimate}
\mynormb{\Delta_h\mathcal{N}^n}
\leq&\, c_0\mynormb{\Delta_h\tilde{u}^{n}}
+2(K_0+1)c_0\mynormb{\Delta_h\tilde{u}^{n-1}}\nonumber\\
\leq &\, c_0\mynormb{\tilde{w}^{n}}+2(K_0+1)c_0\mynormb{\tilde{w}^{n-1}}
+(2K_0+3)c_0c_1h^2\,,
\end{align}
where the equation \eqref{errorSystem-auxiliary2} and the estimate \eqref{globalSpaceError-1}
are used.

Now, taking the inner product of \eqref{errorSystem-auxiliary1OPERATED} by $\tilde{w}_h^{n}$,
one gets
\begin{align}\label{ErrorEstimate-energy}
\myinnerb{\brat{\ffd{\alpha}\tilde{w}}^n,\tilde{w}^{n}}
\leq\myinnerb{\Delta_h\mathcal{N}^n,\tilde{w}^{n}}
+\myinnerb{\brat{\ffd{\alpha}\eta}^n-\Delta_h\xi^n,\tilde{w}^{n}},
\end{align}
because the zero-valued boundary condition in \eqref{errorSystem-auxiliary2OPERATED} leads to
$\myinnerb{\Delta_h\tilde{w}^{n},\tilde{w}^{n}}\leq0$.
Lemma \ref{lemma:Coefficient2-decreasingFastL1} (I) says that the kernels $A^{(n)}_{n-k}$ are decreasing,
so the Cauchy-Schwarz inequality gives
\begin{align*}
\myinnerb{\brat{\ffd{\alpha}\tilde{w}}^n,\tilde{w}^{n}}
\geq&\,A_{0}^{(n)}\mynormt{\tilde{w}^n}^2
-\sum_{k=1}^{n-1}\brab{A_{n-k-1}^{(n)}-A_{n-k}^{(n)}}\mynormt{\tilde{w}^k}\mynormt{\tilde{w}^n}
-A_{n-1}^{(n)}\mynormt{\tilde{w}^0}\mynormt{\tilde{w}^n}\\
=&\,\mynormt{\tilde{w}^n}\kbraB{A_{0}^{(n)}\mynormt{\tilde{w}^n}
-\sum_{k=1}^{n-1}\brab{A_{n-k-1}^{(n)}-A_{n-k}^{(n)}}\mynormt{\tilde{w}^k}
-A_{n-1}^{(n)}\mynormt{\tilde{w}^0}}\\
=&\,\mynormt{\tilde{w}^n}\sum_{k=1}^{n}A_{n-k}^{(n)}\diff\brab{\mynormt{\tilde{w}^k}}\,.
\end{align*}
Thus with the help of Cauchy-Schwarz inequality and \eqref{nonlinearTerm-estimate},
it follows from \eqref{ErrorEstimate-energy} that
\begin{align*}
\sum_{k=1}^{n}A_{n-k}^{(n)}&\,\diff\brab{\mynormt{\tilde{w}^k}}\leq
\mynormb{\Delta_h\mathcal{N}^n}
+\mynormb{\brat{\ffd{\alpha}\eta}^n-\Delta_h\xi^n}\\
\leq&\,
c_0\mynormb{\tilde{w}^{n}}+2(K_0+1)c_0\mynormb{\tilde{w}^{n-1}}
+\mynormb{\brat{\ffd{\alpha}\eta}^n-\Delta_h\xi^n}+(2K_0+3)c_0c_1h^2\,.
\end{align*}
Setting the maximum time-step $\tau\leq\tau_0^{*}=1/\sqrt[\alpha]{3\max\{1,\rhomax\}\Gamma(2-\alpha)(2K_0+3)c_0}$,
we apply Theorem \ref{thm: Gronwall} with
$\xi_1^n=\mynormb{\brat{\ffd{\alpha}\eta}^n-\Delta_h\xi^n}$ and $\xi_2^n=(2K_0+3)c_0c_1h^2$  to get
\begin{align*}
\mynormb{\tilde{w}^{n}}\leq&\, E_n
\bra{2\mynormb{\eta^{0}}
+2\max_{1\leq j\leq n}\sum_{k=1}^jP_{j-k}^{(j)}\mynormb{\brat{\ffd{\alpha}\eta}^k-\Delta_h\xi^k}+3(2K_0+3)c_0c_1\omega_{1+\alpha}(t_n)h^2}\\
\leq&\,E_n\braB{\frac{2c_3}{\sigma(1-\alpha)}\tau^{\min\{2-\alpha,\gamma\sigma\}}
+\frac{2c_4}{\sigma}t_{n}^{\alpha}\hat{t}_{n-1}^{\,2}\epsilon}\\
&\,+E_n\braB{2c_1+\frac{2c_2}{\sigma}\hat{t}_n^{\,2}+3(2K_0+3)c_0c_1\omega_{1+\alpha}(t_n)}h^2\leq E_n\mathcal{T}^n,
\end{align*}
where the initial data \eqref{errorSystem-auxiliary2OPERATED} and
the three estimates \eqref{globalSpaceError-1}-\eqref{globalTimeError-gradedGrid} are used.
Then the error equation \eqref{errorSystem-auxiliary2} with \eqref{globalSpaceError-1} imply
that the claimed error estimate \eqref{ErrorEstimate} holds for $k=n$,
\begin{align*}
\mynormt{\Delta_h\tilde{u}^n}
\leq E_n\mathcal{T}^n+c_1h^2\,.
\end{align*}
The principle of induction and the third inequality in \eqref{discreteNormsEmbedding}
give the following result.

\begin{theorem}\label{th:Convergence-nonuniformFastL1Scheme}
Assume that the solution of nonlinear subdiffusion problem \eqref{prob}
with the nonlinear function $f\in C^4(\mathbb{R})$
fulfills the regularity assumption \eqref{SolutionRegularityAssumption}
with $\sigma\in(0,1)\cup(1,2)$.
If the SOE approximation error $\epsilon\leq \epsilon_0$
and the maximum step size $\tau\leq \min\{\tau_0,\tau_1,\tau_0^{*}\}$,
the discrete solution of two-level linearized fast scheme
\eqref{diffScheme},
on the nonuniform time mesh
satisfying \Ass{3} and \Ass{G}, is unconditionally convergent,
\begin{align*}
\mynormb{U^k-u^k}_{\infty}\leq \frac{C_u}{\sigma(1-\alpha)}\max\{1,\rho\}
\bra{\tau^{\min\{2-\alpha,\gamma\sigma\}}+h^2+\epsilon}, \quad 1\leq k\leq N.
\end{align*}
It achieves an optimal time accuracy of order $O(\tau^{2-\alpha})$
if $\gamma\geq\max\{1,(2-\alpha)/\sigma\}$.
\end{theorem}

\section{Numerical experiments}\label{sec:numerical}

Two numerical examples are reported here to support our theoretical analysis.
The two-level linearized scheme \eqref{diffScheme}
runs for solving the fractional Fisher equation
$$\fd{\alpha}u=\Delta u+u(1-u)+g(\vecx ,t),\quad (\vecx,t)\in(0,\pi)^2\times(0,T],$$
subject to zero-valued boundary data, with two different initial data and exterior forces:
\begin{itemize}
  \item  (Example 1) $u^0(\vecx )=\sin x \sin y$ and $g(\vecx ,t)=0$
  such that no exact solution is available;
  \item  (Example 2) $g(\vecx ,t)$ is specified such that
  $u(\vecx ,t)=\omega_{\sigma}(t)\sin x \sin y$, $0<\sigma<2$.
\end{itemize}
Note that, Example 2 with the regularity parameter $\sigma$ is set to examine
the sharpness of predicted time accuracy
on nonuniform meshes. Actually, our present theory also fits for the semi\/linear problem
with nonzero force $g(\vecx ,t)\in C(\bar\Omega\times[0,T])$.

\begin{figure}[!ht]
  \centering
  \includegraphics[width=2.9in]{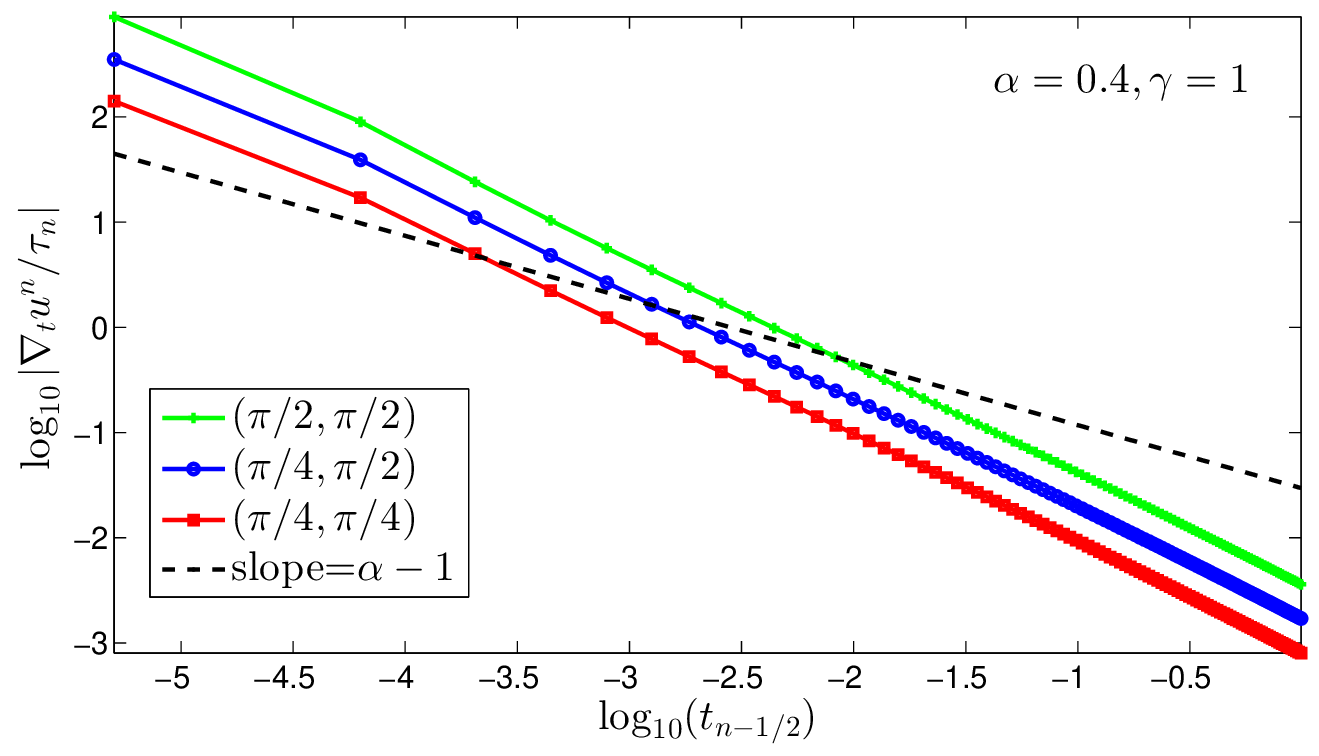}\hfill
  \includegraphics[width=2.9in]{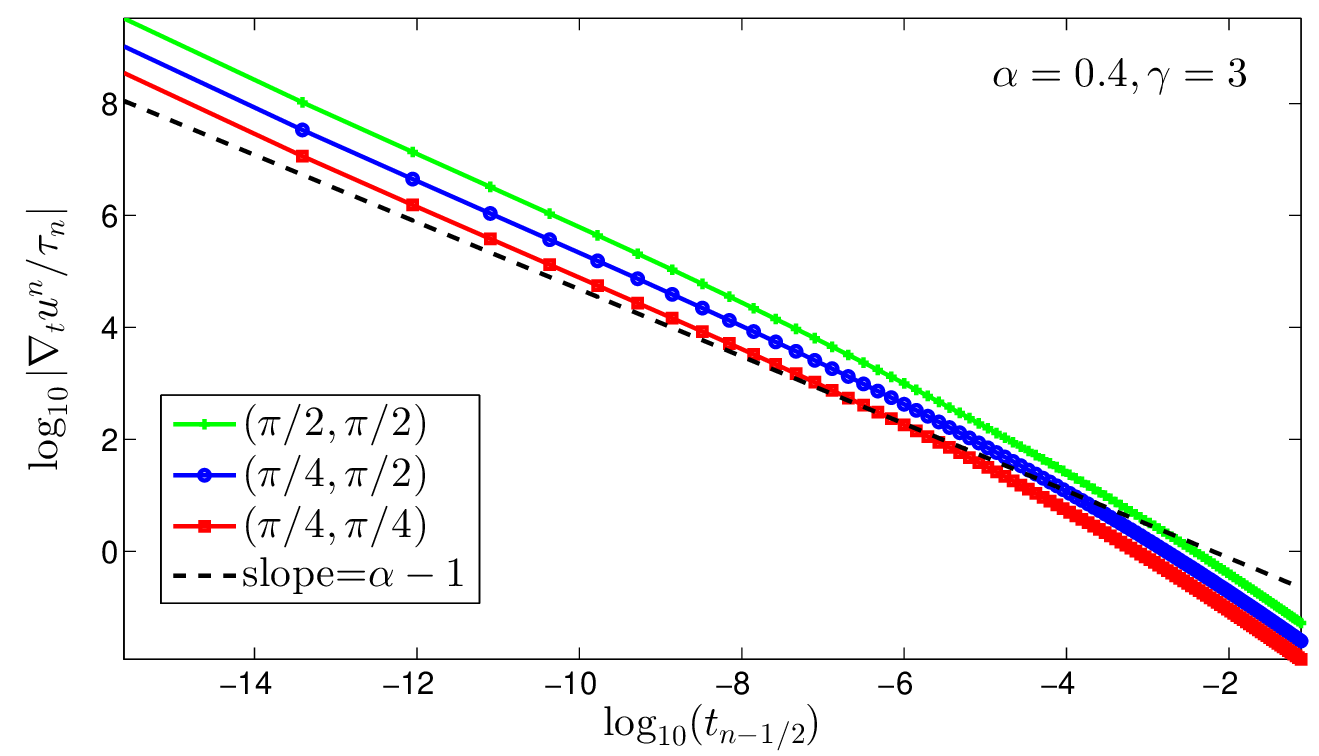}
  \caption{The log-log  plot of difference quotient $\nabla_{\tau}u_h^n/\tau_n$ versus the time
  for Example 1 ($\alpha=0.4$) with two grading parameters $\gamma=1$ (left) and $\gamma=3$ (right).}
  \label{fig:initialsingularityAlpha04}
\end{figure}

\begin{figure}[!ht]
  \centering
  \includegraphics[width=2.9in]{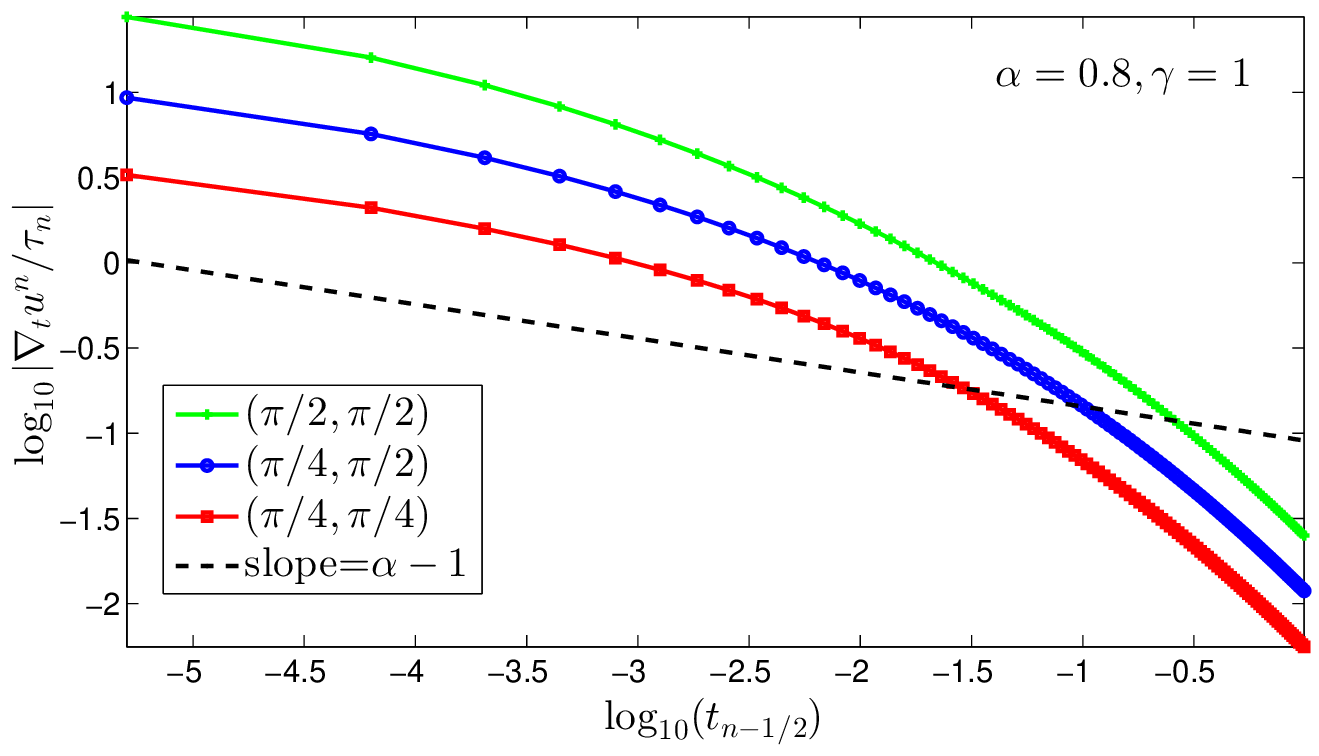}\hfill
  \includegraphics[width=2.9in]{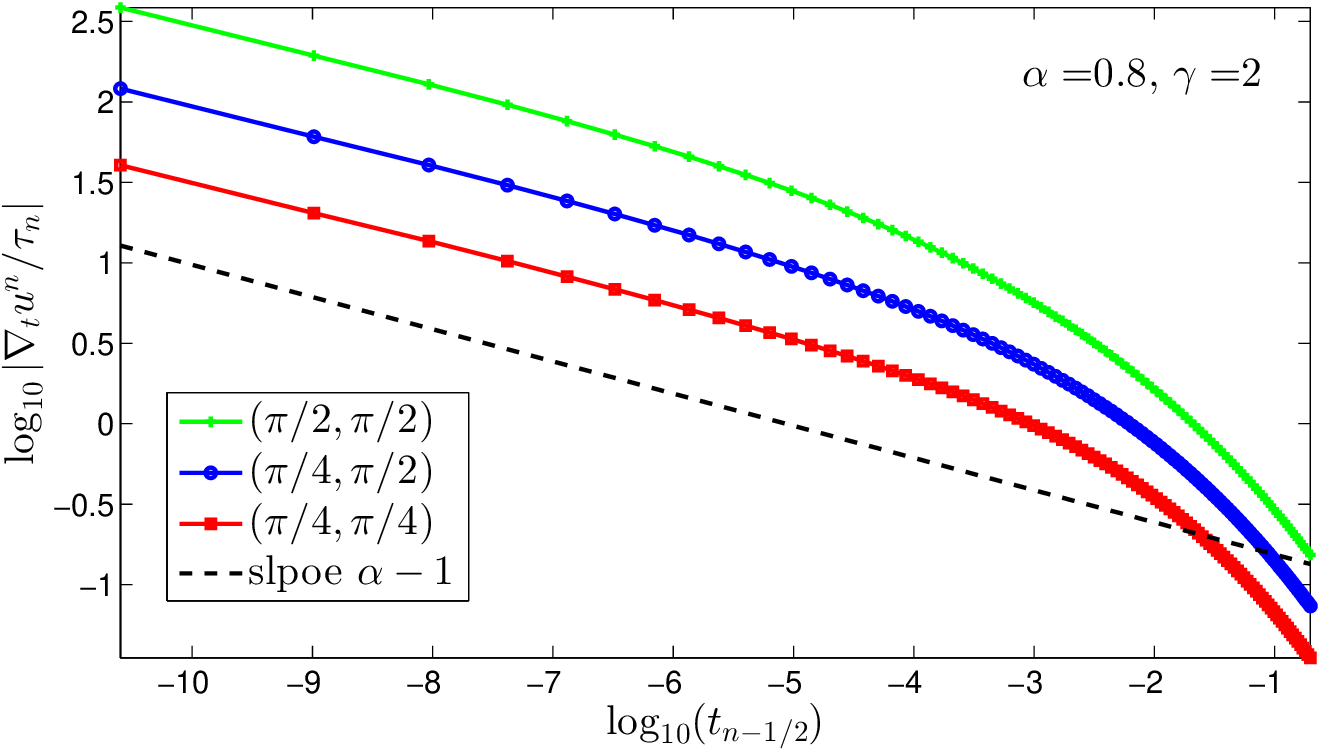}
  \caption{The log-log plot of difference quotient $\nabla_{\tau}u_h^n/\tau_n$ versus the time
  for Example 1 ($\alpha=0.8$) with two grading parameters $\gamma=1$ (left) and $\gamma=2$ (right).}
  \label{fig:initialsingularityAlpha08}
\end{figure}

In our simulations, the spatial domain $\Omega$ is divided uniformly into $M$ parts in each direction $(M_1 = M_2=M)$ and
the time interval $[0,T]$ is divided into two parts $[0, T_0]$ and $[T_0, T]$ with total $N_T$ subintervals.
According to the suggestion in \cite{LiaoLiZhangZhao:2017}, the graded mesh $t_k=T_0\bra{k/N}^{\gamma}$
is applied in the cell $[0, T_0]$ and the uniform mesh with time step size $\tau\geq\tau_{N}$
is used over the remainder interval.
Given certain final time $T$ and a proper number $N_T$,
here we would take $T_0=\min\{1/\gamma,T\}$, $N=\big\lceil\frac{N_T}{T+1-\gamma^{-1}}\big\rceil$ such that
$\tau=\frac{T-T_0}{N_T-N}\geq\frac{T+1-\gamma^{-1}}{N_T}\geq N^{-1}\geq \tau_N.$
Always, the absolute tolerance error of SOE approximation is set to $\epsilon=10^{-12}$
such that the two-level L1 formula \eqref{fastL1formulaNonuniform} is comparable
with the L1 formula \eqref{L1formulaNonuniform} in time accuracy.

In Example 1, we investigate the asymptotic behavior of solution near $t=0$ and the computational efficiency
of the linearized method \eqref{diffScheme}.
Setting $M = 100$, $T=1/\gamma$ and $N_T= 100$,
Figures \ref{fig:initialsingularityAlpha04}-\ref{fig:initialsingularityAlpha08}
depict, in log-log plot, the numerical behaviors
of first-order difference quotient $\nabla_{\tau}u_h^n/\tau_n$ at
three spatial points near the initial time for different
fractional orders and grading parameters. Observations suggest that
$\log\abs{u_t(\vecx ,t)}\approx C_u(\vecx )+(\alpha-1)\log t $ as $t\rightarrow0$, and
the solution is weakly singular near the initial time. Compared with the uniform grid,
the graded mesh always concentrates much more points in the initial
time layer and provides better resolution for the initial singularity.

\begin{figure}[!ht]
  \centering
  \includegraphics[width=3.8in]{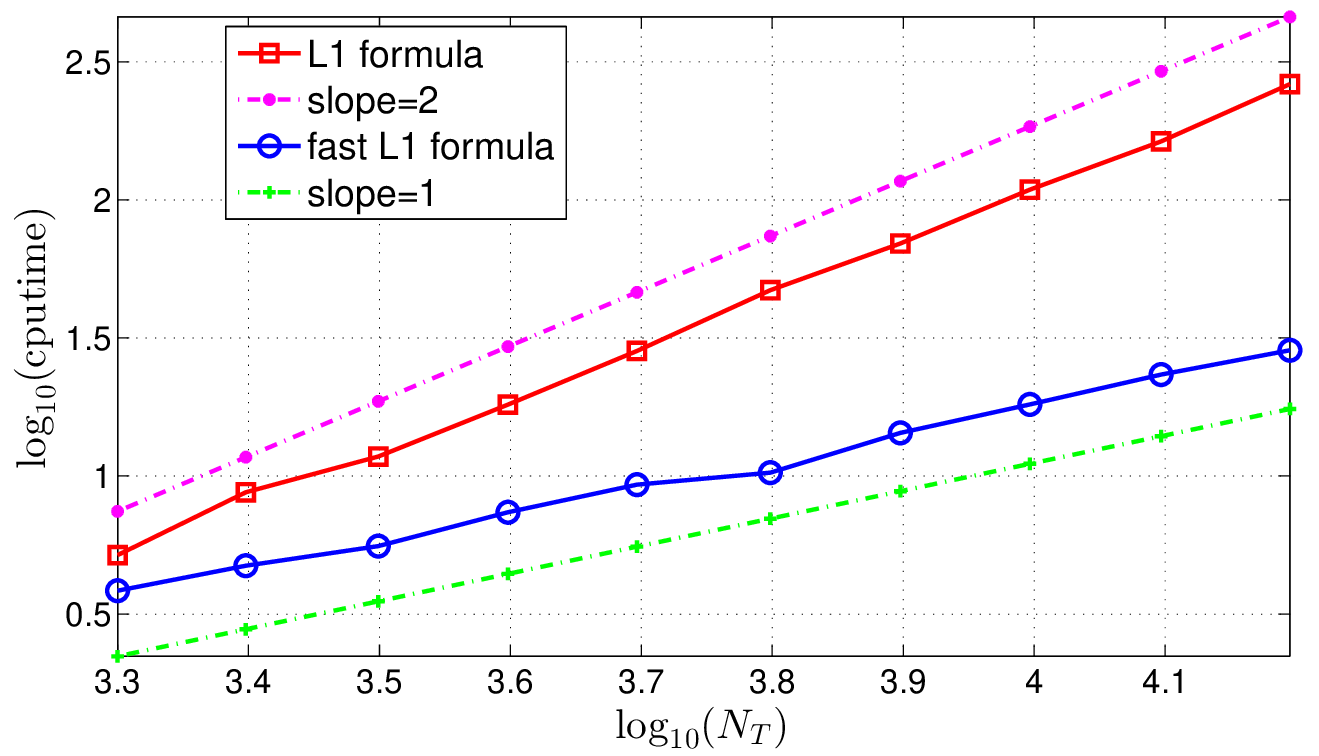}
  \caption{The log-log plot of CPU time versus the total number $N_T$
  of time levels for the linearized method in solving
  Example 1 with two different formulas of Caputo derivative.}
  \label{fig:computationalEfficiency}
\end{figure}

To see the effectiveness of our linearized method \eqref{diffScheme},
we also consider another linearized method by replacing the two-level fast L1 formula $\brat{\ffd{\alpha}u_h}^n$ with
the nonuniform L1 formula  $\brat{\dfd{\alpha}u_h}^n$ defined in \eqref{L1formulaNonuniform}.
Setting $\alpha = 0.5$, $\gamma=2$, and $M = 50$, the two schemes are run for Example 1 to the final time $T = 50$
with different total numbers $N_T$. Figure \ref{fig:computationalEfficiency} shows
the CPU time in seconds for both linearized procedures versus the total number $N_T$ of subintervals.
We observe that the proposed method has almost linear complexity in $N_T$ and is
much faster than the direct scheme using traditional L1 formula.

Since the spatial error $O(h^{2})$ is standard, the time accuracy due to the numerical
approximations of Caputo derivative and nonlinear reaction is examined in Example 2 with $T=1$.
The maximum norm error $e(N,M)=\max_{1\leq l\leq N}\mynormb{U(t_l)-u^l}_{\infty}.$
To test the sharpness of our error estimate,
we consider three different scenarios, respectively, in Tables 5.1-5.3:
\begin{description}
  \item[Table 5.1]: $\sigma=2-\alpha$ and $\gamma=1$ with fractional orders $\alpha=0.4$, $0.6$ and $0.8$.
  \item[Table 5.2]: $\alpha=0.4$ and $\sigma=0.4$ with grid parameters $\gamma=1$,
  $\frac{3}{4}\gamma_{\texttt{opt}}$, $\gamma_{\texttt{opt}}$ and $\frac{5}{4}\gamma_{\texttt{opt}}$.
  \item[Table 5.3]: $\alpha=0.4$ and $\sigma=0.8$ with grid parameters $\gamma=1$,
  $\frac{3}{4}\gamma_{\texttt{opt}}$, $\gamma_{\texttt{opt}}$ and $\frac{5}{4}\gamma_{\texttt{opt}}$.
\end{description}

\begin{table}[!ht]
\begin{center}
\tabcolsep 0pt {Table 5.1 \quad Numerical temporal accuracy for $\sigma=2-\alpha$ and $\gamma=1$} \vspace*{0.5pt}
\def\temptablewidth{1.0\textwidth}
{\rule{\temptablewidth}{1pt}}
\begin{tabular*}{\temptablewidth}{@{\extracolsep{\fill}}ccccccc}
 $N$ &\multicolumn{2}{c}{$\alpha=0.4$,$\sigma=1.6$}&\multicolumn{2}{c}{$\alpha=0.6$,$\sigma=1.4$}
 &\multicolumn{2}{c}{$\alpha=0.8$,$\sigma=1.2$} \\
\cline{2-3}    \cline{4-5} \cline{6-7}& $e(N)$ & Order &$e(N)$ & Order &$e(N)$&Order\\
\hline
50  &5.69e-04 &--  &1.14e-03&--  &2.57e-03&-- \\
100 &1.57e-04 &1.86&4.65e-04&1.30&1.23e-03&1.07\\
200 &4.40e-05 &1.84&1.88e-04&1.31&5.80e-04&1.08\\
400 &1.45e-05 &1.60&7.51e-05&1.32&2.71e-04&1.10\\
800 &5.02e-06 &1.53&2.98e-05&1.34&1.25e-04&1.12\\
\hline
$\min\{\gamma\sigma,2-\alpha\}$& &1.60&&1.40&&1.20 \\
       \end{tabular*}
       {\rule{\temptablewidth}{1pt}}
       \end{center}\label{table:uniform2}
       \end{table}

Tables 5.1 lists the solution errors, for $\sigma=2-\alpha$, on the gradually
refined grids with the coarsest grid of $N=50$. Numerical data indicates that
the optimal time order is of about $O(\tau^{2-\alpha})$,
which dominates the spatial error $O(h^{2})$.
Always, we take $M=N$ in Tables 5.1-5.3 such that $e(N,M)\approx e(N)$.
The experimental rate (listed as Order in tables) of
convergence is estimated by observing that $e(N)\approx C_u\tau^{\beta}$  and then
$\beta\approx\log_{2}\kbra{{e(N)}/{e(2N)}}.$

       \begin{table}[!ht]
\begin{center}
\tabcolsep 0pt {Table 5.2 \quad Numerical temporal accuracy for
$\alpha=0.4$, $\sigma=0.4$ and $\gamma_{\texttt{opt}}=4$} \vspace*{0.5pt}
\def\temptablewidth{1.0\textwidth}
{\rule{\temptablewidth}{1pt}}
\begin{tabular*}{\temptablewidth}{@{\extracolsep{\fill}}ccccccccc}
 $N$ &\multicolumn{2}{c}{$\gamma=1$}&\multicolumn{2}{c}{$\gamma=3$}
 &\multicolumn{2}{c}{$\gamma=4$}&\multicolumn{2}{c}{$\gamma=5$} \\
\cline{2-3}    \cline{4-5} \cline{6-7}\cline{8-9}& $e(N)$ & Order
&$e(N)$ & Order &$e(N)$&Order &$e(N)$&Order\\
\hline
50  &5.47e-02 &--  &3.82e-03&--  &1.65e-03&--  &1.32e-03&--\\
100 &4.64e-02 &0.24&1.68e-03&1.18&5.78e-04&1.52&4.60e-04&1.52 \\
200 &3.78e-02 &0.30&7.36e-04&1.19&1.99e-04&1.54&1.58e-04&1.54 \\
400 &3.00e-02 &0.33&3.21e-04&1.20&6.78e-05&1.55&5.37e-05&1.56 \\
800 &2.34e-02 &0.36&1.40e-04&1.20&2.30e-05&1.56&1.81e-05&1.57\\
\hline
$\min\{\gamma\sigma,2-\alpha\}$& &0.40&&1.20&&1.60&&1.60 \\
       \end{tabular*}
       {\rule{\temptablewidth}{1pt}}
       \end{center}\label{table:smallSigma}
       \end{table}

\begin{table}[!ht]
\begin{center}
\tabcolsep 0pt {Table 5.3 \quad Numerical temporal accuracy for
$\alpha=0.4$, $\sigma=0.8$ and $\gamma_{\texttt{opt}}=2$} \vspace*{0.5pt}
\def\temptablewidth{1.0\textwidth}
{\rule{\temptablewidth}{1pt}}
\begin{tabular*}{\temptablewidth}{@{\extracolsep{\fill}}ccccccccc}
 $N$ &\multicolumn{2}{c}{$\gamma=1$}&\multicolumn{2}{c}{$\gamma=3/2$}
 &\multicolumn{2}{c}{$\gamma=2$}&\multicolumn{2}{c}{$\gamma=5/2$} \\
\cline{2-3}    \cline{4-5} \cline{6-7}\cline{8-9}& $e(N)$ & Order &$e(N)$ & Order &$e(N)$&Order &$e(N)$&Order\\
\hline
50  &3.46e-03 &--  &8.72e-04&--  &5.80e-04&--  &7.52e-04&--\\
100 &2.20e-03 &0.65&3.93e-04&1.15&1.39e-04&2.08&1.77e-04&2.08 \\
200 &1.34e-03 &0.72&1.75e-04&1.17&3.80e-05&1.87&4.06e-05&2.13 \\
400 &7.95e-04 &0.75&7.70e-05&1.18&1.32e-05&1.53&8.88e-06&2.19 \\
600 &5.83e-04 &0.77&4.76e-05&1.19&7.06e-06&1.54&4.22e-06&1.55\\
800 &4.67e-04 &0.77&3.38e-05&1.19&4.52e-06&1.55&2.70e-06&1.55\\
\hline
$\min\{\gamma\sigma,2-\alpha\}$& &0.80&&1.20&&1.60&&1.60 \\
       \end{tabular*}
       {\rule{\temptablewidth}{1pt}}
       \end{center}\label{table:middleSigma}
       \end{table}

Numerical results in Tables 5.2-5.3 (with $\alpha=0.4$ and $\sigma<2-\alpha$)
support the predicted time accuracy
in Theorem \ref{th:Convergence-nonuniformFastL1Scheme} on
the smoothly graded mesh $t_k=T(k/N)^{\gamma}$.
In the case of a uniform mesh $(\gamma=1)$, the solution
is accurate of order $O(\tau^{\sigma})$, and the nonuniform meshes
improve the numerical precision and convergence rate of solution evidently.
The optimal time accuracy $O(\tau^{2-\alpha})$ is observed when the grid parameter $\gamma\geq(2-\alpha)/\sigma$.

\section*{Acknowledgements}
The authors gratefully thank Professor Martin Stynes for his valuable discussions
and fruitful suggestions during the preparation of this paper.
Hong-lin Liao would also thanks for the hospitality of Beijing CSRC during the period of his visit.

\appendix
\section{Proof of Lemma \ref{lemma:BoundednessOf2rdDiscretederivate}}
\def\theequation{A.\arabic{equation}}
\setcounter{equation}{0}

\begin{proof}Consider $F(\psi)=\psi$ firstly. It is easy to check that, at point $\vecx_h=(x_i,y_j)\in\Omega_h$,
\begin{align*}
\delta_x^2(\psi_{ij}v_{ij})=\psi_{ij}\brab{\delta_x^2v_{ij}}+\delta_x\psi_{i-\frac12,j}\brab{\delta_xv_{i-\frac12,j}}
+\delta_x\psi_{i+\frac12,j}\brab{\delta_xv_{i+\frac12,j}}+v_{ij}\brab{\delta_x^2\psi_{ij}}\,,
\end{align*}
so that $\mynormt{\delta_x^2(\psi v)}\leq C_0\bra{\mynormt{v}+\mynormt{\delta_xv}+\mynormt{\delta_x^2v}}$.
Similarly, $\mynormt{\delta_y^2(\psi v)}\leq C_0\bra{\mynormt{v}+\mynormt{\delta_yv}+\mynormt{\delta_y^2v}}.$
Moreover, one has
$\mynormt{\delta_y\delta_x(\psi v)}\leq
C_0\bra{\mynormt{v}+\mynormt{\delta_xv}+\mynormt{\delta_yv}+\mynormt{\delta_y\delta_xv}}$,
due to the fact
\begin{align*}
\delta_y\delta_x(\psi_{i-\frac12,j-\frac12}v_{i-\frac12,j-\frac12})
=&\,\psi_{i-\frac12,j-\frac12}\brab{\delta_y\delta_xv_{i-\frac12,j-\frac12}}
+\delta_y\psi_{i-\frac12,j-\frac12}\brab{\delta_xv_{i-\frac12,j-\frac12}}\\
&\,+\delta_x\psi_{i-\frac12,j-\frac12}\brab{\delta_yv_{i-\frac12,j-\frac12}}
+\brab{\delta_y\delta_x\psi_{i-\frac12,j-\frac12}}v_{i-\frac12,j-\frac12}\,.
\end{align*}
Noticing that
$\mynorm{\Delta_hv}^2=\mynormt{\delta_x^2v}^2+2\mynormt{\delta_x\delta_yv}^2+\mynormt{\delta_y^2v}^2$,
we apply the embedding inequalities in \eqref{discreteNormsEmbedding} to obtain,
also see \cite[Lemma 2.2]{LiaoSunShi:2010},
\begin{align*}
\mynorm{\Delta_h(\psi v)}\leq C_u\bra{\mynorm{v}+\mynorm{\Delta_hv}}\leq C_{F}\mynorm{\Delta_hv},
\end{align*}
where the constant $C_F$ is dependent on $C_0$ and $C_{\Omega}$.
For the general case $F\in C^2(\mathbb{R})$, one has
\begin{align*}
\delta_x^2\kbrab{F(\psi_{ij})v_{ij}}=&\,F(\psi_{ij})\brab{\delta_x^2v_{ij}}
+\delta_xF(\psi_{i-\frac12,j})\brab{\delta_xv_{i-\frac12,j}}\\
&\,+\delta_xF(\psi_{i+\frac12,j})\brab{\delta_xv_{i+\frac12,j}}+v_{ij}\kbrab{\delta_x^2F(\psi_{ij})}\,.
\end{align*}
The formula of Taylor expansion with integral remainder gives
\begin{align*}
\delta_xF(\psi_{i-\frac12,j})=&\,\brab{F(\psi_{ij})-F(\psi_{i-1,j})}/h_1
=\delta_x\psi_{i-\frac12,j}\int_0^1F'\brab{s \psi_{ij}+(1-s)\psi_{i-1,j}}\zd{s}\,,\\
\delta_x^2F(\psi_{ij})=&\,\brab{\delta_x^2\psi_{ij}}F'(\psi_{ij})
+\brab{\delta_x\psi_{i-\frac12,j}}^2\int_0^1F''\brab{s \psi_{ij}+(1-s)\psi_{i-1,j}}(1-s)\zd{s}\\
&\,\hspace{2.4cm}+\brab{\delta_x\psi_{i+\frac12,j}}^2\int_0^1F''\brab{s \psi_{ij}+(1-s)\psi_{i+1,j}}(1-s)\zd{s}\,,
\end{align*}
such that $\mynormt{\delta_xF(\psi)}\leq C_{F}$ and $\mynormt{\delta_x^2F(\psi)}\leq C_{F}\,$.
Therefore, simple calculations arrive at
\begin{align*}
\mynormt{\delta_x^2\kbra{F(\psi) v)}}\leq C_{F}\bra{\mynormt{v}+\mynormt{\delta_xv}+\mynormt{\delta_x^2v}}\,.
\end{align*}
By presenting similar arguments as those in the above simple case,
it is straightforward to get claimed estimate and complete the proof.
\end{proof}

\end{document}